\RequirePackage{fix-cm} 
\documentclass[a4paper, twoside, reqno, dvips, 12pt]{amsart}
\usepackage{fixltx2e}   

\usepackage{etex}

\usepackage[utf8,latin1]{inputenc}
\usepackage[T1]{fontenc}

\usepackage{eucal}
\usepackage{esint}
\usepackage{dsfont}
\usepackage{xspace}
\usepackage{amsgen}
\usepackage{amsthm}
\usepackage{amssymb}
\usepackage{amsmath}
\usepackage{upgreek}
\usepackage{wasysym}
\usepackage{amsfonts}
\usepackage{mathrsfs}
\usepackage{stmaryrd}
\usepackage{mathtools}

\usepackage{a4wide}

\usepackage[dvipsnames, table]{xcolor}
\definecolor{bckg}{RGB}{20.8, 20.8, 20.8}
\definecolor{oneblue}{rgb}{0.0, 0.0, 0.85}
\definecolor{Lightblue}{RGB}{214, 214, 214}
\definecolor{bluepigment}{rgb}{0.2, 0.2, 0.6}
\definecolor{charcoal}{rgb}{0.21, 0.27, 0.31}
\definecolor{denimblue}{rgb}{0.08, 0.38, 0.74}
\definecolor{Lightgray}{rgb}{0.89, 0.89, 0.89}
\definecolor{darkgrey}{rgb}{0.273, 0.281, 0.30}
\definecolor{darkelectricblue}{rgb}{0.33, 0.41, 0.47}

\usepackage[sort&compress, comma, square, numbers]{natbib}

\usepackage[framed, numbered, autolinebreaks, useliterate]{mcode}

\usepackage{psfrag}
\usepackage{graphicx}
\usepackage{subfigure}
\usepackage{morefloats}
\usepackage{indentfirst}

\usepackage{acronym}
\usepackage{microtype}
\usepackage[labelsep=period,%
            labelfont={bf,sf,color=bluepigment},%
            justification=raggedright]{caption}

\usepackage[perpage, symbol]{footmisc}

\usepackage[usenames, dvipsnames, pdf]{pstricks}
\usepackage{epsfig}
\usepackage{pst-grad} 
\usepackage{pst-plot} 

\usepackage{tikz-cd}

\usepackage[colorlinks,
           urlcolor=oneblue,
           linkcolor=denimblue,
           citecolor=NavyBlue,
           bookmarksopen=false,
           pdfpagemode=UseNone,
           pagebackref]{hyperref}

\usepackage[explicit]{titlesec}

\titleformat{\section}
  {\color{NavyBlue}\Large\sffamily\bfseries}
  {}
  {0em}
  {\colorbox{bckg!5}{\parbox{\dimexpr\linewidth-2\fboxsep\relax}{\centering\thesection. #1}}}
  [\vspace*{0.33em}]

\titleformat{name=\section,numberless}
  {\color{NavyBlue}\Large\sffamily\bfseries}
  {}
  {0.0em}
  {\colorbox{bckg!10}{\parbox{\dimexpr\linewidth-2\fboxsep\relax}{\centering#1}}}
  [\vspace*{0.33em}]

\titleformat{\subsection}
  {\color{NavyBlue}\large\sffamily\bfseries}
  {}
  {0.0em}
  {\colorbox{bckg!5}{\parbox{\dimexpr\linewidth-2\fboxsep\relax}{\centering\thesubsection. #1}}}
  [\vspace*{0.33em}]

\titleformat{name=\subsection,numberless}
  {\color{NavyBlue}\Large\sffamily\bfseries}
  {}
  {0em}
  {\colorbox{bckg!10}{\parbox{\dimexpr\linewidth-2\fboxsep\relax}{\centering#1}}}
  [\vspace*{0.33em}]

\titleformat{\subsubsection}
  {\color{bluepigment}\sffamily\normalsize\bfseries}
  {\thesubsubsection}
  {0.5em}
  {#1}
  [\vspace*{0.33em}]

\titleformat{\paragraph}[runin]
  {\color{bluepigment}\sffamily\small\bfseries}
  {}
  {0em}
  {#1}

\titlespacing{\section}{1.0em}{1.5em plus 2pt minus 2pt}%
{1.0em plus 2pt minus 2pt}[0em]
\titlespacing{\subsection}{1.0em}{1.5em plus 2pt minus 2pt}%
{1.0em}[0em]
\titlespacing{\subsubsection}{1.0em}{1.5em plus 2pt minus 2pt}%
{1.0em plus 2pt minus 2pt}[0em]

\usepackage{titletoc}

\contentsmargin{0.5em}
\newlength{\tocsep} 
\titlecontents{section}[\tocsep]
  {\addvspace{10pt}\bfseries\sffamily}
  {\contentslabel[\thecontentslabel]{\tocsep}}
  {}
  {\ \titlerule*[0.75pc]{.}\ \thecontentspage}
  []
\titlecontents{subsection}[\tocsep]
  {\addvspace{8pt}\sffamily}
  {\contentslabel[\thecontentslabel]{\tocsep}}
  {}
  {\ \titlerule*[0.5pc]{.}\ \thecontentspage}
  []
\titlecontents*{subsubsection}[\tocsep]
  {\addvspace{2pt}\footnotesize\sffamily}
  {}
  {}
  {\ \titlerule*[0.35pc]{.}\ \thecontentspage}
  [\\*]

\makeatletter
\def\@setauthors{%
  \begingroup
  \def\thanks{\protect\thanks@warning}%
  \trivlist
  \centering\footnotesize \@topsep30\p@\relax
  \advance\@topsep by -\baselineskip
  \item\relax
  \author@andify\authors
  \def\\{\protect\linebreak}%
  \textsc{\normalsize\textcolor{darkelectricblue}{\authors}}%
  \ifx\@empty\contribs
  \else
    ,\penalty-3 \space \@setcontribs
    \@closetoccontribs
  \fi
  \endtrivlist
  \endgroup
}
\def\@settitle{\begin{center}%
  \baselineskip14\p@\relax
    \bfseries
    \textsc{\Large\textcolor{charcoal}{\@title}}
  \end{center}%
}
\makeatother

\usepackage{enumitem}
\setlist[description]{%
  topsep=30pt,               
  itemsep=5pt,               
  font={\bfseries\sffamily\color{NavyBlue}}, 
}

\usepackage{fancyhdr}
\usepackage{lastpage}

\newcommand*\Title{\textcolor{bluepigment}{A brief introduction to pseudo-spectral methods}}
\newcommand*\Authors{\textcolor{bluepigment}{Denys~Dutykh}}
\newcommand*{\plogo}{\textcolor{gray}{{\texttt{arXiv.org} / \textsc{hal}}}} 

\numberwithin{equation}{section}

\newtheorem{exo}{Exercise}
\newtheorem{remark}{Remark}
\newtheorem{deff}{Definition}
\newtheorem{theorem}{Theorem}
\newtheorem{corollary}{Corollary}

\newcommand{\up}[1]{$^{\mathrm{\small\textsf{#1}}}$} 

\newcommand{\s}{{\sf s}}

\newcommand{\hv}{\hat{v}}
\newcommand{\Th}{\hat{T}}

\newcommand{\E}{\mathbb{E}}

\newcommand{\C}{\mathds{C}}
\newcommand{\N}{\mathds{N}}
\newcommand{\R}{\mathds{R}}
\newcommand{\Z}{\mathds{Z}}

\newcommand{\Id}{\mathbb{I}}
\newcommand{\A}{\mathscr{A}}
\newcommand{\B}{\mathscr{B}}
\newcommand{\D}{\mathcal{D}}
\newcommand{\I}{\mathcal{I}}
\newcommand{\K}{\mathcal{K}}

\newcommand{\ud}{\mathrm{d}}
\newcommand{\ui}{\mathrm{i}}
\newcommand{\ue}{\mathrm{e}}
\newcommand{\F}{\mathcal{F}}
\newcommand{\Tt}{\mathbb{T}}
\newcommand{\U}{\mathcal{U}}
\newcommand{\V}{\mathcal{V}}
\newcommand{\W}{\mathcal{W}}
\newcommand{\Gr}{\mathcal{G}}
\newcommand{\Pl}{\mathcal{P}}

\newcommand{\Ll}{\mathscr{L}}
\newcommand{\Nn}{\mathscr{N}}

\renewcommand{\P}{\mathbb{P}}

\newcommand{\Rr}{\mathscr{R}}
\renewcommand{\beta}{\upbeta}
\newcommand{\Hh}{\mathscr{H}}

\renewcommand{\leq}{\leqslant}
\renewcommand{\geq}{\geqslant}
\newcommand{\eps}{\varepsilon}
\renewcommand{\O}{\mathcal{O}}
\renewcommand{\L}{\mathcal{L}}
\renewcommand{\H}{\mathcal{H}}
\renewcommand{\S}{\mathcal{S}}

\renewcommand{\alpha}{\upalpha}

\newcommand{\n}{\boldsymbol{n}}
\newcommand{\x}{\boldsymbol{x}}
\newcommand{\vO}{\boldsymbol{0}}

\renewcommand{\k}{\boldsymbol{k}}

\newcommand{\thetah}{\hat{\theta}}
\newcommand{\const}{\mathrm{const}}
\newcommand{\phis}{\tilde{\varphi}}
\newcommand{\Sh}{\hat{\mathcal{S}}}

\newcommand{\Wb}{\boldsymbol{\mathcal{W}}}

\newcommand{\ie}{\emph{i.e.}}
\newcommand{\eg}{\emph{e.g.}~}
\newcommand{\etc}{\emph{etc.}~}

\renewcommand{\sim}{\thicksim}
\newcommand{\sech}{\mathrm{sech}}
\newcommand{\scal}{\boldsymbol{\cdot}}
\newcommand{\grad}{\boldsymbol{\nabla}}
\newcommand{\abs}[1]{\lvert\, #1\, \rvert}

\newcommand{\norm}[1]{\lVert\, #1\, \rVert}

\newcommand{\pd}[2]{\dfrac{\partial #1}{\partial\/ #2}}
\newcommand{\eqdef}{\mathop{\stackrel{\mathrm{def}}{:=}}}

\newcommand{\od}[2]{\dfrac{\mathrm{d} #1}{\mathrm{d}\/#2}}

\newcommand{\half}{{\textstyle{1\over2}}}

\usepackage{acronym}
\acrodef{bvp}[BVP]{Boundary Value Problem}
\acrodef{NSWE}{Nonlinear Shallow Water Equations}

\begin{document}

\title[\Title]{A brief introduction to pseudo-spectral methods: application to diffusion problems}

\author[D.~Dutykh]{Denys Dutykh$^*$}
\address{LAMA, UMR 5127 CNRS, Universit\'e Savoie Mont Blanc, Campus Scientifique, 
73376 Le Bourget-du-Lac Cedex, France}
\email{Denys.Dutykh@univ-savoie.fr}
\urladdr{http://www.denys-dutykh.com/}
\thanks{$^*$ Corresponding author}

\keywords{pseudo-spectral methods; approximation theory; diffusion; parabolic equations}


\begin{titlepage}
\thispagestyle{empty} 
\noindent
{\Large Denys \textsc{Dutykh}}\\
{\it\textcolor{gray}{CNRS--LAMA, Universit\'e Savoie Mont Blanc, France}}
\\[0.16\textheight]

\vspace*{0.99cm}

\colorbox{Lightblue}{
  \parbox[t]{1.0\textwidth}{
    \centering\huge\sc
    \vspace*{0.7cm}
    
    \textcolor{bluepigment}{A brief introduction to pseudo-spectral methods: application to diffusion problems}
    
    \vspace*{0.7cm}
  }
}

\vfill 

\raggedleft     
{\large \plogo} 
\end{titlepage}


\newpage
\thispagestyle{empty} 
\par\vspace*{\fill}   
\begin{flushright} 
{\textcolor{denimblue}{\textsc{Last modified:}} \today}
\end{flushright}


\newpage
\maketitle
\thispagestyle{empty}


\begin{abstract}

The topic of these notes could be easily expanded into a full one-semester course. Nevertheless we shall try to give some flavour along with theoretical bases of spectral and pseudo-spectral methods. The main focus is made on \textsc{Fourier}-type discretizations, even if some indications on how to handle non-periodic problems via \textsc{Tchebyshev} and \textsc{Legendre} approaches are made as well. The applications presented here are diffusion-type problems in accordance with the topics of the PhD school.

\bigskip
\noindent \textbf{\keywordsname:} pseudo-spectral methods; approximation theory; diffusion; parabolic equations \\

\smallskip
\noindent \textbf{MSC:} \subjclass[2010]{ 65M70, 65N35 (primary), 80M22, 76M22 (secondary)}
\smallskip \\
\noindent \textbf{PACS:} \subjclass[2010]{ 47.11.Kb (primary), 44.35.+c (secondary)}

\end{abstract}


\newpage
\tableofcontents
\thispagestyle{empty}


\newpage

\section{General introduction}

\bigskip
\begin{flushright}{\slshape
    A numerical simulation is like sex. \\
    If it is good, then it is great. \\
    If it is bad, then it is still better than nothing.} \\ \medskip
    --- Dr.~\textsc{D} (paraphrasing Dr.~\textsc{Z})
\end{flushright}
\bigskip
\begin{flushright}{\slshape
    People think they do not understand Mathematics, \\
    but it is all about how you explain it to them.} \\ \medskip
    --- I. M.~\textsc{Gelfand}
\end{flushright}
\bigskip\bigskip

This document represents a collection of brief notes of lectures delivered by the Author at the International PhD school ``\emph{Numerical Methods for Diffusion Phenomena in Building Physics: Theory and Practice}'' which took place in Pontifical Catholic University of Parana (PUCPR, Curitiba, Brazil) in April, 2016. The Author of the present document is grateful to the Organizing Committee of this school (in particular to Professors Nathan~\textsc{Mendes} and Marx~\textsc{Chhay}) for giving him an opportunity to lecture there. The present document should be considered as a supplementary material to the Lectures delivered by the Author at this PhD school. The exposition below is partially based on \cite{Fornberg1996} and some other references mentioned in the manuscript.

This lecture is organized as follows. First, we present some theoretical bases behind spectral discretizations in Section~\ref{sec:intro}. An application to a problem stemming from the building physics is given in Section~\ref{sec:heat}. Finally, we give some indications for the further reading in Section~\ref{sec:indic}. This document contains also a certain number of Appendices directly or indirectly related to spectral methods. For instance, in Appendix~\ref{app:tcheb} we give some useful identities about \textsc{Tchebyshev} polynomials and in Appendix~\ref{app:trefftz} we give some flavour of \textsc{Trefftz} methods, which are essentially forgotten nowadays. We decided to include also some history of the diffusion in Sciences in general and in Physics in particular. You will find this information in Appendix~\ref{app:hist}. Finally, we prepared also an Appendix~\ref{app:mc} devoted to the Monte--Carlo methods to simulate numerically diffusion processes. The Author admits that it is not related to spectral methods, but nevertheless we decided to include it since these methods remain essentially unknown in the community of numerical methods for Partial Differential Equations (PDEs).


\section{Introduction to spectral methods}
\label{sec:intro}

Consider for simplicity a 1D compact\footnote{A domain $\U \subseteq \R^d$, $d \geq 1$ is compact if it is bounded and closed. For a more general definition of compactness we refer to any course in General Topology.} domain, \eg $\;\U\ =\ [-1,\, 1]$ an Ordinary or Partial Differential Equation (ODE or PDE) on it
\begin{equation}\label{eq:pde}
  \Ll u\ =\ g\,, \qquad x\ \in\ \U
\end{equation}
where $u(x,\,t)$ (or just $u(x)$ in the ODE case) is the solution\footnote{For simplicity we consider in this Section scalar equations only. The generalizations for systems is straightforward, since one can apply the same discretization for every individual component of the solution.} which satisfies some additional (boundary) conditions at $x\ =\ \pm 1$ depending on the (linear or nonlinear) operator $\Ll\ =\ \Ll(\partial_t,\, \partial_x,\, \partial_{xx},\, \ldots)$. Function $g(x,\,t)$ is a source term which does not depend on the solution $u$. For example, if \eqref{eq:pde} is the classical heat equation in the free space (\ie~without heat sources), then
\begin{equation*}
  \Ll\ \equiv\ \partial_t\ -\ \nu\,\partial_{xx}\,, \qquad
  g\ \equiv\ 0\,,
\end{equation*}
where $\nu \in \R$ is the diffusion coefficient. If the solution $u(x)$ is steady (\ie~time-independent), we deal with an ODE. In any case, in the present document we focus only on the discretization in space. The time discretization is discussed deeper in the lecture of Dr.~Marx~\textsc{Chhay}. A brief reminder of some most useful numerical techniques for ODEs is given in Appendix~\ref{app:odes}.

The idea behind a spectral method is to approximate a solution $u(x,\,t)$ by a finite sum
\begin{equation}\label{eq:expa}
  u(x,\,t)\ \approx\ u_n(x,\,t)\ =\ \sum_{k\, =\, 0}^{n}\, v_k\,(t)\, \phi_k\,(x),
\end{equation}
where $\{\phi_k(x)\}_{k=0}^{\infty}$ is the set of basis functions. In ODE case all $v_k(t)\ \equiv\ \const$. The main question which arises is how to choose the basis functions? Once the choice of $\{\phi_k(x)\}_{k\,=\,0}^{\infty}$ is made, the second question appears: how to determine the expansion coefficients $v_k(t)$?

Here we shall implicitly assume that the function $u(x,\,t)$ is \emph{smooth}. Only in this case the full potential of spectral methods can be exploited\footnote{We have to say that pseudo-spectral methods can be applied to problems featuring eventually discontinuous solutions as well (\eg hyperbolic conservation laws). However, it is out of scope of the present lecture devoted rather to parabolic problems. As a side remark, the Author would add that the advantage of pseudo-spectral methods for nonlinear hyperbolic equations is not so clear comparing to modern high-resolution shock-capturing schemes \cite{Leer2006}.}. The concept of \emph{smoothness} in Mathematics is ambiguous since there are various classes of smooth functions:
\begin{equation*}
  C^{\,p}(\U)\ \supseteq \ C^{\,\infty}(\U)\ \supseteq\ \A^{\,\infty}(\D), \qquad (p\ \in\ \Z^+, \quad \U\ \subseteq\ \D\ \subseteq\ \C).
\end{equation*}
The last sequence of inclusions needs perhaps some explanations:
\begin{description}
  \item[$C^{\,p}(\U)$] the class of functions $f: \U\ \mapsto\ \R$ having at least $p\ \geq\ 1$ continuous derivatives.
  \item[$C^{\,\infty}(\U)$] the class of functions $f: \U\ \mapsto\ \R$ having infinitely many continuous derivatives.
  \item[$\A^{\,\infty}(\D)$] the class of functions $f: \D\ \mapsto\ \C$ analytical (holomorphic) in a domain $\D\ \subseteq\ \C$ containing the segment $\U$ in its interior.
\end{description}
For example, the \textsc{Runge}\footnote{Carl David Tolm\'e \textsc{Runge} (1856 -- 1927) is a German Mathematician and Physicist. A PhD student of Karl~\textsc{Weierstra\ss}.} function $f(x)\ =\ \dfrac{1}{1 + 25\, x^2}$ is in $C^{\,\infty}\bigl([-1,\,1]\bigr)$, but not analytical in the complex plain $\A^{\,\infty}(\C)$.

\subsection{Choice of the basis}
\label{sec:basis}

A successful expansion basis meets the following requirements:

\begin{enumerate}
  \item~[\textcolor{NavyBlue}{\textbf{Convergence}}] The approximations $u_n(x,\,t)$ should converge rapidly to $u(x,\,t)$ as $n \to \infty$
  \item~[\textcolor{NavyBlue}{\textbf{Differentiation}}] Given coefficients $\{v_k(t)\}_{k=0}^{n}$, it should be easy to determine another set of coefficients\footnote{Here the prime does not mean a time derivative!} $\{v_k'(t)\}_{k=0}^{n}$ such that
  \begin{equation*}
    \pd{u_n}{x}\ =\ \sum_{k\, =\, 0}^{n}\, v_k(t)\, \od{\phi_k(x)}{x}\ \leadsto\ \sum_{k\, =\, 0}^{n}\, v_k^{\,\prime}(t)\, \phi_k(x)\,.
  \end{equation*}
  \item~[\textcolor{NavyBlue}{\textbf{Transformation}}] The computation of expansion coefficients $\{v_k\}_{k=0}^{n}$ from function values $\bigl\{u(x_i,t)\bigr\}_{i=0}^{n}$ and the reconstruction of solution values in nodes from the set of coefficients $\{v_k\}_{k=0}^{n}$ should be easy, \ie~the conversion between two data sets is algorithmically efficient
  \begin{equation*}
    \bigl\{u(x_i,t)\bigr\}_{i=0}^{n} \quad \leftrightarrows \quad \{v_k\}_{k=0}^{n}\,.
  \end{equation*}
\end{enumerate}

\subsubsection{Periodic problems}

For periodic problems it is straightforward to propose a basis which satisfies the requirements (1)--(3) above. It consists of \emph{trigonometric polynomials}:
\begin{equation}\label{eq:trig}
  u_n(x, t)\ =\ a_0(t)\ +\ \sum_{k\, =\, 1}^{n}\,\bigl\{a_k(t)\,\cos(k \pi x)\ +\ b_k(t)\,\sin(k \pi x)\bigr\}.
\end{equation}
The first two points are explained in elementary courses of analysis, while the requirement (3) is possible thanks to the invention of the Fast \textsc{Fourier} Transform (FFT) algorithm first by \textsc{Gau}\ss{} and later by \textsc{Cooley} \& \textsc{Tukey} in 1965 \cite{Cooley1965}.

\begin{remark}
The use of trigonometric bases such as \eqref{eq:trig} or $\{\ue^{\ui\, k \pi x}\}_{k\in\Z}$ for heat conduction problems has been initiated by J.~\textsc{Fourier}\footnote{Jean-Baptiste Joseph \textsc{Fourier} (1768 -- 1830) is a French Physicist and Mathematician. In particular he accompanied Napoleon \textsc{Bonaparte} on his campaign to Egypt as a scientific adviser.} (1822) \cite{Fourier1822} even if \textsc{Fourier} series have been known well before \textsc{Fourier}. There is the so-called \textsc{Arnold}'s\footnote{Vladimir~\textsc{Arnold} (1937 -- 2010), a prominent Soviet/Russian mathematician. Please, read his books!} principle which states that in Mathematics nothing is named after its true inventor. The question the Author would like to rise is whether \textsc{Arnold}'s principle is applicable to it-self?
\end{remark}

\begin{remark}
The term `Spectral methods' can be now explained. It comes from the fact that the solution $u(x,\,t)$ is expanded into a series of orthogonal eigenfunctions of some linear operator $\L$ (with partial or ordinary derivatives). In this way, the numerical solution is related to its spectrum, thus justifying the name `spectral methods'. For example, if we take the \textsc{Laplace} operator $\L\ =\ -\grad^2\ \equiv\ -\sum_{j\, =\, 1}^{d} \pd{^2}{x_j^2}$ on the periodic domain $[\,0,\,2\,\pi)^{\,d}$, its spectrum consists of \textsc{Fourier} modes:
\begin{equation*}
  -\grad^2 \ue^{-\ui\,\k\scal\x}\ =\ \abs{\k}^2\,\ue^{-\ui\,\k\scal\x}\,.
\end{equation*}
In this way we obtain naturally the \textsc{Fourier} analysis and \textsc{Fourier}-type pseudo-spectral methods.
\end{remark}

\subsubsection{Non-periodic problems}

\begin{figure}
  \centering
  \includegraphics[width=0.99\textwidth]{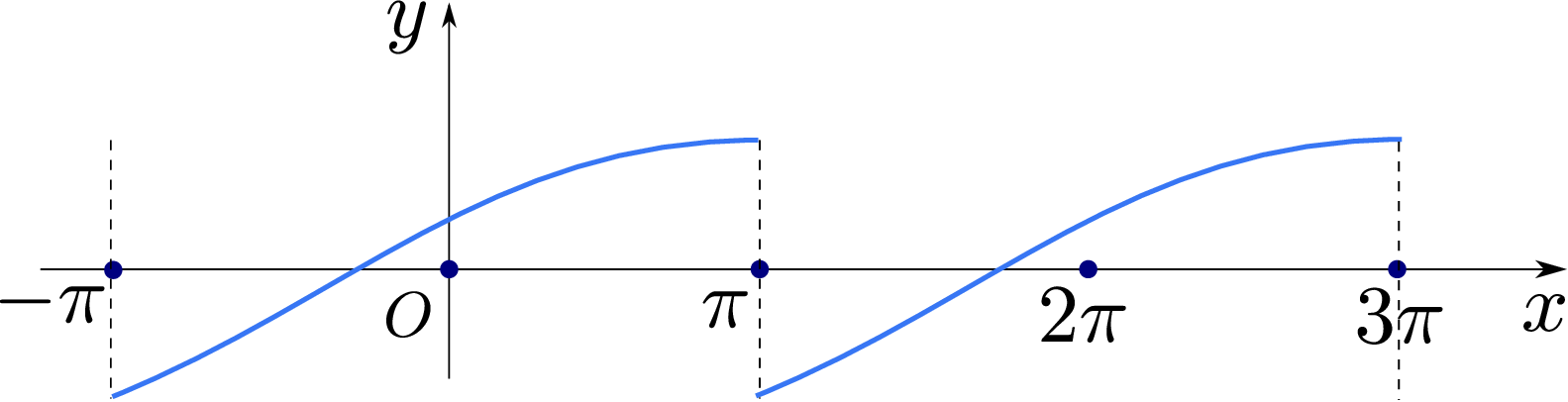}
  \caption{\small\em Periodisation of a smooth continuous function defined on $[-\pi,\,\pi]$.}
  \label{fig:period}
\end{figure}

The trigonometric basis \eqref{eq:trig} fails to work for general non-periodic problems essentially because of the failure of requirement (1). Indeed, the artificial discontinuities arising after the periodisation (see Figure~\ref{fig:period}) make the \textsc{Fourier} coefficients $v_n$ decay as $\O(n^{-1})$ when $n \to \infty$.

The analytic series (or \textsc{Taylor}-like expansions) represent another interesting alternative:
\begin{equation}\label{eq:taylor}
  u_n(x,\,t)\ =\ \sum_{k\, =\, 1}^{n}\, v_k(t)\, x^k\,.
\end{equation}
The problem is that \textsc{Taylor}-type expansions \eqref{eq:taylor} work well only for a very limited class of functions. Namely, they satisfy the requirement (1) only for functions analytical\footnote{The analyticity property is understood here in the sense of the Complex Analysis.} in the unit disc $\D_1(\vO)\ \subseteq\ \C$. For instance, the celebrated \textsc{Runge}'s function $f(x)\ =\ \dfrac{1}{1 + 25 x^2}$ fails to satisfy this condition because of two imaginary poles located at $z\ =\ \pm \dfrac{\ui}{5}$.

The successful examples of polynomial bases are given by \emph{orthogonal polynomials}, which arise naturally in various contexts:

\begin{description}
  
  \item[Numerical integration] Optimal numerical integration formulas achieve a high accuracy by using zeros of orthogonal polynomials as nodes.

  \item[Sturm--Liouville problem] \textsc{Jacobi}\footnote{Carl Gustav Jacob~\textsc{Jacobi} (1804 -- 1851) is a German Mathematician.} polynomials arise as eigenfunctions to singular \textsc{Sturm}\footnote{Jacques Charles Fran\c{c}ois \textsc{Sturm} (1803 -- 1855) is a French Mathematician who was born in Geneva which was a part of France at that time.}--\textsc{Liouville}\footnote{Joseph \textsc{Liouville} (1809 -- 1882) is a French Mathematician who founded the \emph{Journal de Math\'matiques Pures et Appliqu\'ees}.} problems.
  
  \item[Approximation in $L_2$] Truncated expansions in \textsc{Legendre}\footnote{Adrien-Marie~\textsc{Legendre} (1752 -- 1833) is a French Mathematician. Not to be confused with an obscure politician Louis \textsc{Legendre}. Their portraits are often confused.} polynomials are optimal approximants in the $L_2$-norm.

  \item[Approximation in $L_\infty$] Truncated expansions in \textsc{Tchebyshev}\footnote{Pafnuty~\textsc{Tchebyshev} (1821 -- 1894), a Russian mathematician who contributed to many fields of Mathematics from number theory to probabilities and numerical analysis.} polynomials are optimal approximants in the $L_\infty$-norm.

\end{description}

\begin{figure}
  \centering
  \includegraphics[width=0.75\textwidth]{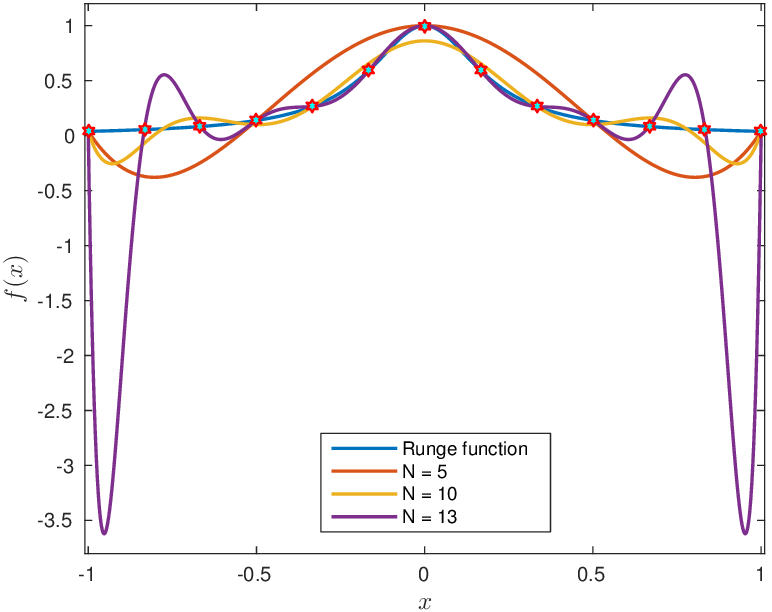}
  \caption{\small\em Interpolation of the \textsc{Runge} function $f(x)\ =\ \dfrac{1}{1 + 25 x^2}$ on a sequence of successively refined of uniform grids. One can observe the divergence phenomenon close to the interval boundaries $x\ =\ \pm 1$.}
  \label{fig:runge}
\end{figure}

\begin{figure}
  \centering
  \includegraphics[width=0.75\textwidth]{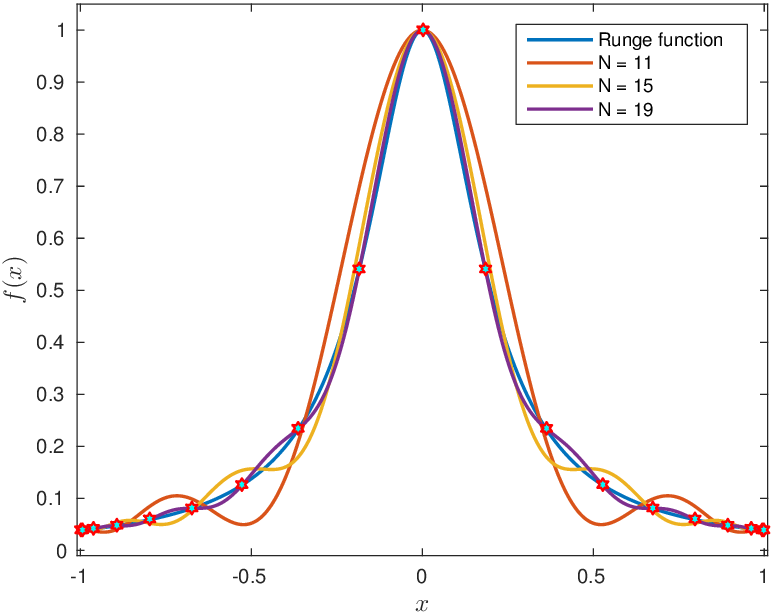}
  \caption{\small\em Interpolation of the \textsc{Runge} function $f(x)\ =\ \dfrac{1}{1 + 25 x^2}$ on a sequence of successively refined of \textsc{Tchebyshev}'s grids. One can see that the oscillations present in Figure~\ref{fig:runge} disappear and the interpolant seems to converge to the interpolated function $f(x)\,$.}
  \label{fig:runge2}
\end{figure}

Each of the topics above deserves a separate course to be covered. Here we content just to provide this information as facts, which can be deepened later, if necessary. We would just like to quote \textsc{Bernardi} \& \textsc{Maday} (1991):
\begin{quote}
  \it
  We do think that the corner stone of collocation techniques is the choice of the collocation nodes [$\,$\dots] in spectral methods these are always built from the nodes of a Gau\ss~quadrature formula.
\end{quote}
Consequently, extrema (and zeros) of \textsc{Tchebyshev} (and some other orthogonal) polynomials play a very important r\^ole in the Numerical Analysis (NA). \textsc{Tchebyshev} nodes are given explicitly by
\begin{equation}\label{eq:nodes}
  x_k\ =\ -\cos\,\Bigl(\frac{\pi k}{N}\Bigr)\ \in\ [-1,\,1]\,, \qquad k\ =\ 0,\, 1,\, 2,\, \ldots,\, N\,.
\end{equation}
There is a result saying that using nodes \eqref{eq:nodes} as interpolation points gives an interpolant which is not too far from the optimal polynomial $\Pl^{\mathrm{Opt}}_N$, \ie
\begin{equation*}
  \norm{f\ -\ \Pl^{\mathrm{Tch}}_N}_{\infty}\ \leq\ \bigl(1\ +\ \Lambda^{\mathrm{Tch}}_N\bigr)\; \norm{f\ -\ \Pl^{\mathrm{Opt}}_N}_{\infty},
\end{equation*}
where $f:\ [-1,\,1]\ \mapsto\ \R$ is the function that we interpolate and $\Pl^{\mathrm{Tch}}_N$ is the interpolation polynomial constructed on nodes \eqref{eq:nodes}. Here $\Lambda^{\mathrm{Tch}}_N$ is the so-called \textsc{Lebesgue}\footnote{Henri~\textsc{Lebesgue} (1875 -- 1941), a French mathematician most known for the theory of integration having its name.} constant for \textsc{Tchebyshev} interpolation polynomials of degree $N$. Notice that $\Lambda^{\mathrm{Tch}}_N$ does not depend on the function $f$ being interpolated. The \textsc{Lebesgue} constants for various interpolation techniques have the following asymptotic behaviour
\begin{equation*}
  \Lambda^{\mathrm{Tch}}_N\ \sim\ \O(\log N)\,, \qquad
  \Lambda^{\mathrm{Leg}}_N\ \sim\ \O(\sqrt{N})\,, \qquad
  \Lambda^{\mathrm{Uni}}_N\ \sim\ \O\Bigl(\frac{2^N}{N\log N}\Bigr)\,,
\end{equation*}
where $\Lambda^{\mathrm{Leg}}_N$ and $\Lambda^{\mathrm{Uni}}_N$ are \textsc{Lebesgue} constants for \textsc{Legendre} and uniform nodes distributions correspondingly. In particular one can see that the uniform node distribution is simply disastrous. This phenomenon is another `avatar' of the so-called \textsc{Runge} phenomenon illustrated in Figure~\ref{fig:runge}. The performance of \textsc{Tchebyshev}'s nodes is shown in Figure~\ref{fig:runge2}.

\begin{remark}
As it was shown by \textsc{V\'ertesi} (1990) \cite{Vertesi1990}, the \textsc{Lebesgue} constant for \textsc{Tchebyshev} distribution of nodes is very close for the smallest possible \textsc{Lebesgue} constant:
\begin{align*}
  \Lambda^{\mathrm{Tch}}_N\ &=\ \frac{2}{\pi}\;\Bigl(\,\ln N\ +\ \gamma\ +\ \ln\frac{8}{\pi}\,\Bigr)\ +\ o\,(1), \\
  \Lambda^{\mathrm{min}}_N\ &=\ \frac{2}{\pi}\;\Bigl(\,\ln N\ +\ \gamma\ +\ \ln\frac{4}{\pi}\,\Bigr)\ +\ o\,(1)\,.
\end{align*}
where $\gamma\ \approx\ 0.5772156649\dots$ is the \textsc{Euler}\footnote{Leonhard~\textsc{Euler} (1707 -- 1783) is a great mathematician who was born in Switzerland and worked all his life in Saint-Petersburg.}--\textsc{Mascheroni}\footnote{Lorenzo~\textsc{Mascheroni} (1750 -- 1800) was an Italian mathematician who worked in Pavia.} constant.
\end{remark}

\bigskip
\paragraph*{The Lebesgue constant.}

It is worth to explain better the important notion of the \textsc{Lebesgue} constant and how it appears in the theory of interpolation (below we follow \cite{Solin2005}). Let us take an arbitrary (valid) nodes distribution $\{x_i\}_{i\,=\,0}^{n}$ in domain $\U$, \ie
\begin{equation*}
  \forall i\ =\ 1,\,2,\,\ldots\,,\,n\,:\ x_i\ \in\ \U\,, \qquad \forall j\ \neq\ i\,:\ x_i\ \neq\ x_j\,.
\end{equation*}
For any continuous function $u\ \in\ C(\U)$ there exists a unique interpolating polynomial $\Pl(x)\ \in\ \P_{n}\,[\,\R\,]$ of degree $n\ =\ \deg \Pl\,$:
\begin{equation*}
  \Pl_n(x_i)\ =\ u_i\ \equiv\ u\,(x_i)\,, \qquad i\ =\ 1,\,2,\,\ldots\,,\,n\,.
\end{equation*}
This interpolating polynomial $\Pl_n(x)$ will be denoted as $\Id_{\,n}\,[\,u\,]\,$. We introduce also the best possible approximating polynomial $\Pl^{\,\ast}\ \in\ \P_{n}\,[\,\R\,]\,$:
\begin{equation*}
  \norm{u\ -\ \Pl^{\,\ast}}_{\infty}\ =\ \inf_{\Pl\ \in\ \P_{n}\,[\,\R\,]}\,\norm{u\ -\ \Pl}_{\infty}\,.
\end{equation*}
It is not obliged that $\Pl^{\,\ast}\ \equiv\ \Id_{\,n}\,[\,u\,]\,$. However, since $\Pl^{\,\ast}\ \in\ \P_{n}\,[\,\R\,]$ we necessarily have $\Pl^{\,\ast}\ \equiv\ \Id_{\,n}\,[\,u\,]\,$. Therefore, the following inequalities hold:
\begin{multline*}
  \norm{u\ -\ \Id_{\,n}\,[\,u\,]}_{\infty}\ =\ \norm{u\ -\ \Pl^{\,\ast}\ +\ \Pl^{\,\ast}\ -\ \Id_{\,n}\,[\,u\,]}_{\infty}\ \equiv\\ \norm{u\ -\ \Pl^{\,\ast}\ +\ \Id_{\,n}\,[\,\Pl^{\,\ast}\,]\ -\ \Id_{\,n}\,[\,u\,]}_{\infty}\ \leq\ \norm{u\ -\ \Pl^{\,\ast}}_{\infty}\ +\ \norm{\Id_{\,n}\,[\,\Pl^{\,\ast}\,]\ -\ \Id_{\,n}\,[\,u\,]}_{\infty}\\
  \leq\ \norm{u\ -\ \Pl^{\,\ast}}_{\infty}\ +\ \norm{\Id_{\,n}}\cdot\norm{\Pl^{\,\ast}\ -\ u}_{\infty}\ \leq\ \bigl(1\ +\ \norm{\Id_{\,n}}\bigr)\cdot\norm{u\ -\ \Pl^{\,\ast}}_{\infty}\,.
\end{multline*}
The norm of a linear operator $\Id_{\,n}\,[\cdot]$ is defined as
\begin{equation*}
  \norm{\Id_{\,n}}\ \eqdef\ \sup_{\norm{u}_{\infty}\ =\ 1}\,\norm{\Id_{\,n}\,[\,u\,]}_{\infty}\,.
\end{equation*}
The aforementioned norm of the interpolation operator $\Id_{\,n}\,[\cdot]$ is called the \textsc{Lebesgue} constant for the set of nodes $\{x_i\}_{i\,=\,0}^{n}\,$. In multiple dimensions the \textsc{Lebesgue} constant depends also on the shape of domain $\U$ additionally to the nodes distribution $\{x_i\}_{i\,=\,0}^{n}\,$. As we mentioned above, obvious choices such as the uniform distribution of nodes is disastrous, since it yields the exponential growth of the \textsc{Lebesgue} constant $\Lambda^{\mathrm{Uni}}_n\ \sim\ \O\Bigl(\frac{2^n}{n\log n}\Bigr)\,$, where $n$ is the degree of the interpolating polynomial.

The estimation of the \textsc{Lebesgue} constant in multi-dimensional non-\textsc{Cartesian} domains is a problem essentially open nowadays. The most precious information is to find the nodes distribution, for example over a triangle, which minimizes the \textsc{Lebesgue} constant. This knowledge would be crucial for the design of new spectral elements \cite{Solin2005}. The locations of nodes which minimize the magnitude of the \textsc{Lebesgue} constant $\Lambda_{n}\ \to\ \min$ for a given polynomial space $\P_{n}\,[\,\R\,]$ are called \textsc{Lebesgue} points. Nowadays, the best known points on quadrangles are tensor products of \textsc{Gau\ss}--\textsc{Lobatto}\footnote{Rehuel~\textsc{Lobatto} is a Dutch Mathematician born in a Portuguese family.} or \textsc{Tchebyshev} nodes. On triangles the \textsc{Fekete}\footnote{Michael~\textsc{Fekete} (1886 -- 1957) is a Hungarian Mathematician who did his PhD under the supervision of Lip\'ot~\textsc{Fej\'er}. M.~\textsc{Fekete} gave also some private tutoring to J\'anos\textsc{Neumann} known today as John~\textsc{von Neumann}.} points seem currently to be the best choice.

\bigskip
\paragraph*{Intermediate conclusions.}

The remarks above show that \textsc{Tchebyshev} polynomials satisfy the requirement (1) for a useful spectral basis. The requirement (2) is met due to derivative recursion formulas, which can be easily demonstrated (see also Appendix~\ref{app:tcheb}):
\begin{equation*}
  T_k^{\,\prime}(x)\ =\ 2\,k\,\sum_{n\,=\,0}^{\left[\frac{k-1}{2}\right]}\;\frac{1}{\delta_{k-1-2\,n}}\;T_{k-1-2\,n}\,(x)\,, \qquad
  \delta_k\ =\ \begin{dcases}
    \ 2\,,&\ k\ =\ 0\,, \\
    \ 1\,,&\ k\ \neq\ 0\,.
  \end{dcases}
\end{equation*}
Finally, the requirement (3) is satisfied as well thanks to the Fast Cosine \textsc{Fourier} Transform (FCFT) (a variant of FFT). It allows to compute spectral coefficients $\{v_k\}_{k\, =\, 0}^{N}$ from the node values $\bigl\{u_n(x_i)\bigr\}_{i\, =\, 0}^{N}$ and vice versa. Consequently, \textsc{Tchebyshev} polynomials have become an almost universal choice for non-periodic problems. These methods in 1D have been implemented in the \textsc{Matlab} toolbox \texttt{Chebfun}\footnote{\texttt{Chebfun}'s team is lead by Nick~\textsc{Trefethen}.} (and \texttt{Chebfun2} in 2D).

\begin{remark}
To Author's knowledge, \textsc{Legendre} polynomials found some applications in the construction of Spectral Element Method (SEM) bases \cite{Solin2005}.
\end{remark}

\begin{remark}
We have to mention that high-order polynomial interpolants have a bad reputation in the Numerical Analysis (NA) community for the two following reasons:
\begin{enumerate}
  \item Because of the \textsc{Runge} phenomenon
  \item and because of the following
  \begin{theorem}
  For any node density distribution function there exists a continuous function such that the $L_\infty$-norm of the interpolation error tends to infinity when the number of nodes $N\ \to\ +\infty$.
  \end{theorem}
\end{enumerate}
However, this bad reputation is not justified. For instance, the \textsc{Runge} phenomenon is completely suppressed by \textsc{Tchebyshev} nodes distribution. The low or moderate \textsc{Lebesgue} constant $\Lambda_N$ ensures that we are not too far from the optimal polynomial.
\end{remark}

\bigskip
\paragraph*{Infinite domains.}

The infinite domains can be handled, first of all, by periodisation in the context of \textsc{Fourier}-type methods as explained above. However, truly infinite domains can be handled by various approaches:
\begin{itemize}
  \item \textsc{Hermite} polynomials
  \item Sinc functions
  \item (almost) Rational functions
  \item Change of variables, \eg 
  \begin{equation*}
    x:\ [-\pi,\, \pi]\ \mapsto\ \R\,, \qquad x(q)\ =\ \ell\;\tan\,\Bigl(\,\frac{q}{2}\,\Bigr), \qquad \ell\ \in\ \R^+
  \end{equation*}
  \item Truncation of a domain $(-\infty\,, \infty)\ \rightsquigarrow\ [-\ell,\, \ell]$ followed by a change of variables such as
  \begin{equation*}
    x:\ [-\pi,\, \pi]\ \mapsto\ [-\ell,\,\ell]\,, \qquad x(q)\ =\ \frac{\ell\, q}{\pi}\,.
  \end{equation*}
\end{itemize}
Two last items are followed by a conventional pseudo-spectral discretization on a finite domain $[-\pi,\, \pi]$.

A \textsc{Hermite}\footnote{Charles \textsc{Hermite} (1822 -- 1901), a French mathematician whose ``\emph{Cours d'analyse}'' represent a lot of interest even today. He was also the PhD advisor of another great French mathematician --- Henri~\textsc{Poincar\'e}.}-type pseudo-spectral method can be briefly sketched as follows. First of all, we construct recursively the family of \textsc{Hermite}'s polynomials:
\begin{equation*}
  \H_0(x)\ =\ 1\,, \qquad \H_1(x)\ =\ 2\,x, \qquad
  \H_{n+1}(x)\ =\ 2\, x\, \H_n(x)\ -\ 2\, n\, \H_{n-1}(x), \qquad
  n\ \geq\ 1\,.
\end{equation*}
Then, we construct \textsc{Hermite}'s functions, which form an \emph{orthonormal} basis on $\R$:
\begin{equation*}
  \Hh_n(x)\ \eqdef\ \frac{1}{\sqrt{2^n \cdot n!}}\; \H_n(x)\; \ue^{-\frac{x^2}{2}}.
\end{equation*}
Derivatives of expansions in \textsc{Hermite} functions can be easily re-expanded using the relation:
\begin{equation*}
  \Hh_n^{\,\prime}(x)\ =\ -\;\sqrt{\frac{n+1}{2}}\;\Hh_{n+1}(x)\ +\ \sqrt{\frac{n}{2}}\;\Hh_{n-1}(x)\,.
\end{equation*}
However, the lack of a Fast \textsc{Hermite} Transform (FHT) implies the use of \emph{differentiation matrices} in numerical implementations. There are also some concerns about a poor convergence rate when the number of modes $N$ is increased. The infinite domains remain very challenging \emph{in silico}.

\subsection{Determining expansion coefficients}
\label{sec:coeff}

There are three main techniques to find the spectral expansion coefficients $\{v_k\}_{k = 0}^{N}$. In order to explain them, let us introduce first the residual function on a trial solution $u_n(x,\,t)$:
\begin{equation*}
  \Rr[u_n](x,\, t)\ \eqdef\ \Bigl[\Ll u_n\ -\ g\Bigr]\,(x,\, t).
\end{equation*}
Typically, we shall evaluate the residual $\Rr$ on the expansion \eqref{eq:expa} and the residual norm $\norm{\Rr}$ is generally considered as a measure of the approximate solution quality. The goal is to keep the residual as small as possible across the domain $\U$ (in our particular case $\U\ =\ [-1,\, 1]$).

So, we can mention here at least five approaches to determine the expansion coefficients:
\begin{description}
  \item[Tau--Lanczos] Spectral coefficients $\{v_k\}_{k = 0}^{N}$ are selected such that the boundary conditions are satisfied identically and the residual $\Rr[u_n]$ is orthogonal to as many basis functions $\phi_k(x)$ as possible.
\end{description}
\begin{description}
  \item[Galerkin] First the basis functions are recombined $\bigl\{\phi_k(x)\bigr\}_{k = 0}^{N}\ \rightsquigarrow \bigl\{\phis_k(x)\bigr\}_{k = 0}^{N}$ so that the boundary conditions are satisfied identically. Then, the coefficients $\{v_k\}_{k = 0}^{N}$ are found so that the residual $\Rr[u_n]$ be orthogonal to as many of \emph{new} basis functions $\bigl\{\tilde{\phi}_k(x)\bigr\}_{k = 0}^{N}$ as possible.
\end{description}
\begin{description}
  \item[Collocation] This approach is similar to the Tau--Lanczos method concerning the boundary conditions: spectral coefficients $\{v_k\}_{k = 0}^{N}$ are selected such that the boundary conditions are satisfied. The rest of coefficients is determined so that the residual $\Rr[u_n](x,\,t)$ vanishes at as many (thoroughly chosen) spatial locations as possible.
\end{description}
\begin{description}
  \item[Petrov--Galerkin] It is a variant of Galerkin method in which the residual $\Rr[u_n]$ is made orthogonal to a set of functions, which is different from the approximation space basis $\bigl\{\phi_k(x)\bigr\}_{k = 0}^{N}$.
\end{description}
\begin{description}
  \item[Least squares] Various least square-type approaches are used when for some reason the number of coefficients to be determined is different from the number of conditions which can be imposed. Below we shall avoid such pathological situations.
\end{description}

The Tau--Lanczos technique was proposed by C.~\textsc{Lanczos}\footnote{Cornelius~\textsc{Lanczos} (1893 -- 1974), a Hungarian/American numerical analyst. His books are very pedagogical as well.} in 1938. What we call the Galerkin\footnote{Boris~\textsc{Galerkin} (more precisely romanized as \textsc{Galyorkin}) (1871 -- 1945), a Russian then Soviet Civil Engineer. The Author suggests to read his biography which is comparable to James~\textsc{Bond} (007) movies.} technique was proposed independently first by I.~\textsc{Bubnov}\footnote{Ivan~\textsc{Bubnov} (1872 -- 1919), a Russian naval engineer.} and by W.~\textsc{Ritz}\footnote{Walther~\textsc{Ritz} (1878 -- 1909), a talented Swiss Physicist, who died young in G\"ottingen, Germany.} (one more example of the \textsc{Arnold} principle in action!). So, to respect the historical time line, the method should be called \textsc{Bubnov}--\textsc{Ritz}--\textsc{Galerkin}. Today it is the basis of the Finite Element Method (FEM). The collocation technique was called the \emph{pseudo-spectral method} presumably for the first time by S.~\textsc{Orszag}\footnote{Steven~\textsc{Orszag} (1943 -- 2011), an American numerical analyst, one of the first users of pseudo-spectral methods.} in 1972. The \textsc{Petrov}--\textsc{Galerkin} method was proposed by G.I.~\textsc{Petrov}\footnote{Georgii Ivanovich~\textsc{Petrov} (1912 -- 1987), a Russian fluid mechanician.}. It is used up to now in some convection-dominated problems.

\begin{remark}
The collocation method can be recast into the \textsc{Petrov}--\textsc{Galerkin} framework when we make the residual $\Rr[u_n]$ is made orthogonal to \textsc{Dirac}\footnote{Paul Adrien Maurice \textsc{Dirac} (1902 -- 1984) is a British Theoretical Physicist who predicted theoretically using \textsc{Dirac}'s equation the existence of the positron. One of the founders of the Quantum Mechanics.} singular measures $\bigl\{\delta(x\ -\ x_k)\bigr\}_{k = 1}^{n}$, where $\{x_k\}_{k = 1}^{n}$ are the collocation points.
\end{remark}

For linear problems all these methods work equally well. However, for nonlinear ones (and in the presence of variable coefficients) the pseudo-spectral (collocation) approach is particularly easy to apply since it involves the products of numbers (solution/variable coefficient's values in collocation points) instead of products of expansions, which are much more difficult to handle.

The convergence of pseudo-spectral approximations for very smooth functions is always geometrical, \ie~$\ \sim\ \O(q^{\,N})$, where $N$ is the number of modes. This statement is true for any derivative with the same convergence factor $0\ <\ q\ <\ 1$. However, the periodic pseudo-spectral method converges always faster than its non-periodic counterpart. This conclusion follows from convergence properties of \textsc{Fourier} and \textsc{Tchebyshev} series in the complex domain.

The relative resolution ability of various pseudo-spectral methods can be also quantified in terms of the number of points per wavelength needed to resolve a signal. Indeed, this description is more suitable for wave propagation problems. For periodic \textsc{Fourier}-type methods one needs $2$ points per wavelength. For \textsc{Tchebyshev}-type methods one needs about $\pi$ points. Finally, this number goes up to $6$ nodes per wavelength for uniform grids. However, due to the huge \textsc{Lebesgue} constant $\Lambda_N^{\mathrm{Uni}}$ the uniform grids in pseudo-spectral setting are not usable in practice.


\section{Aliasing, interpolation and truncation}

\begin{figure}
  \centering
  \includegraphics[width=0.99\textwidth]{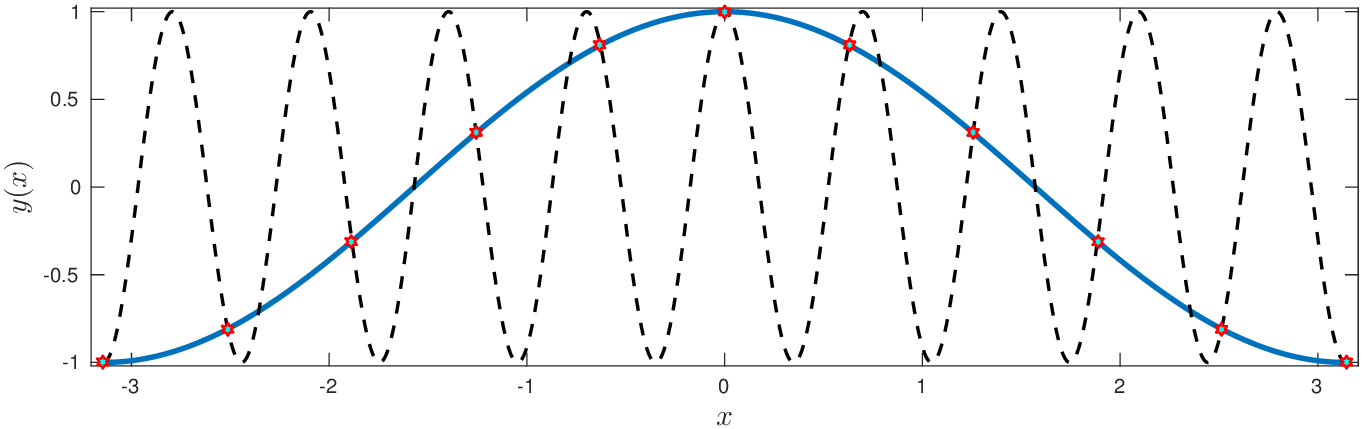}
  \caption{\small\em Illustration of the aliasing phenomenon: two \textsc{Fourier} modes are indistinguishable on the discrete grid. The modes represented here are $\cos(x)$ and $\cos(9x)$ and the discrete grid is composed of $N = 11$ equispaced points on the segment $[-\pi,\,\pi]$.}
  \label{fig:alias}
\end{figure}

Let us take a continuous (and possibly a smooth) function $u(x)$ defined on the interval $\I\ =\ (-\pi,\,\pi)$ and develop it in a \textsc{Fourier} series. In general it will contain the whole spectrum (\ie~the infinite number) of frequencies:
\begin{equation*}
  u(x)\ =\ \sum_{k\, =\, -\infty}^{+\infty} v_k\; \ue^{\ui k x}.
\end{equation*}
Now let us discretize the interval $\I$ with $N$ equispaced collocation points (as we do it in \textsc{Fourier}-type pseudo-spectral methods. In the following we assume $N$ to be odd, \ie~$N\ =\ 2m + 1$. On this discrete grid all modes $\bigl\{\ue^{\ui (k + jN)x}\bigr\}_{j \in \Z}$ are indistinguishable. See Figure~\ref{fig:alias} for an illustration of this phenomenon.

The interpolating trigonometric polynomial on a given grid can be written as
\begin{equation*}
  \Id_N[u]\ =\ \sum_{k\, =\, -m}^{m} \hv_k\;\ue^{\ui k x}\,.
\end{equation*}
Each discrete \textsc{Fourier} coefficient incorporates the contributions of all modes which looks the same on the considered grid:
\begin{equation*}
  \hv_k\ =\ \sum_{j\, =\, -\infty}^{+\infty} v_{k + jN}\,.
\end{equation*}
Let us recall that the polynomial $\Id_N[u]$ takes the prescribed values $\bigl\{u(x_k)\bigr\}_{k\, =\, -m}^{m}$ in the points of the grid $\{x_k\}_{k\, =\, -m}^{m}$. This object is fundamentally different from the truncated \textsc{Fourier} series:
\begin{equation*}
  \Tt_N[u]\ =\ \sum_{k\, =\, -m}^{m} v_k\; \ue^{\ui k x}\,.
\end{equation*}
The difference between these two quantities is known as the \emph{aliasing error}:
\begin{equation*}
  \Rr_N[u]\ \eqdef\ \Id_N[u]\ -\ \Tt_N[u]\ =\ \sum_{k\, =\, -m}^m \sum_{\forall j \neq 0} v_{k + jN}\; \ue^{\ui k x}\,.
\end{equation*}
After applying the Pythagoras theorem\footnote{We can apply the \textsc{Pythagoras} theorem since the aliasing error $\Rr_N[u]$ contains the \textsc{Fourier} modes with numbers $\abs{k} \leq m$, while the reminder $u\ -\ \Tt_N[u]$ contains only the modes with $\abs{k} > m$. Thus, they are orthogonal.}, we obtain
\begin{equation*}
  \norm{u\ -\ \Id_N[u]}^2_{L_2}\ =\ \norm{u\ -\ \Tt_N[u]}^2_{L_2}\ +\ \norm{\Rr_N[u]}^2_{L_2}.
\end{equation*}
Thus, the interpolation error is always larger than the truncation error in the standard $L_2$ norm. The amount of this difference is precisely equal to the committed \emph{aliasing error}. However, we prefer to use in pseudo-spectral methods the interpolation technique because of the Discrete \textsc{Fourier} Transform (DFT), which allows to transform quickly (thanks to the FFT algorithm) from the set of function values in grid points to the set of its interpolation coefficients. So, it is easier to apply an FFT instead of computing $N$ integrals to determine the \textsc{Fourier} series coefficients.

\bigskip
\paragraph*{Nonlinearities.} Let us take the simplest possible nonlinearity --- the product of two functions $u(x)$ and $v(x)$ defined by their truncated \textsc{Fourier} series containing the modes up to $m\,$:
\begin{equation*}
  u(x)\ =\ \sum_{k\, =\, -m}^{m} u_k\; \ue^{\ui k x}\,, \qquad
  v(x)\ =\ \sum_{k\, =\, -m}^{m} v_k\; \ue^{\ui k x}\,.
\end{equation*}
The product of these two functions $w(x)$ can be obtained my multiplying the \textsc{Fourier} series:
\begin{equation*}
  w(x)\ =\ u(x)\cdot v(x)\ =\ \biggl(\,\sum_{k\, =\, -m}^{m} u_k\; \ue^{\ui k x}\,\biggr)\cdot\biggl(\,\sum_{k\, =\, -m}^{m} v_k\; \ue^{\ui k x}\,\biggr)\ \equiv\ \sum_{k\, =\, -2\,m}^{2\,m} w_k\; \ue^{\ui k x}\,.
\end{equation*}
It can be clearly seen that the product contains high order harmonics up to $\ue^{\pm\ui m x}$ which cannot be represented on the initial grid. Thus, they will contribute to the aliasing error explained above.

The aliasing of a nonlinear product can be ingeniously avoided by adopting the so-called $3/2$\up{th} rule whose \textsc{Matlab} implementation is given below. This function assumes that input vectors are \textsc{Fourier} coefficients of functions $u(x)$ and $v(x)\,$. The resulting vector contains (anti-aliased) \textsc{Fourier} coefficients of their product $w(x)\ =\ u(x)\cdot v(x)\,$.
\begin{lstlisting}
function w_hat = AntiAlias(u_hat, v_hat)
  N         = length(u_hat);
  M         = 3*N/2; % 3/2th rule
  u_hat_pad = [u_hat(1:N/2) zeros(1, M-N) u_hat(N/2+1:end)];
  v_hat_pad = [v_hat(1:N/2) zeros(1, M-N) v_hat(N/2+1:end)];
  u_pad     = ifft(u_hat_pad);
  v_pad     = ifft(v_hat_pad);
  w_pad     = u_pad.*v_pad;
  w_pad_hat = fft(w_pad);
  w_hat     = 3/2*[w_pad_hat(1:N/2) w_pad_hat(M-N/2+1:M)];
end % AntiAlias()
\end{lstlisting}
The main idea behind is to complete vectors of \textsc{Fourier} coefficients by a sufficient number of zeros (\ie~the so-called zero padding technique) so that in the physical space the product $u(x)\cdot v(x)$ can be fully resolved. The final step consists in extracting $m$ relevant \textsc{Fourier} coefficients \cite{Trefethen2000, Uecker2009}.

\begin{remark}
To Author's knowledge, the development of efficient and rigorously justified anti-aliasing rules for other types of nonlinearities such as the division, square root, \etc is an open problem.
\end{remark}


\subsection{Example of a second order boundary value problem}

Consider the following second order Boundary Value Problem (BVP) on the interval $\I\ =\ [-1,\, 1]$:
\begin{equation}\label{eq:bvp}
  u_{xx}\ +\ u_x\ -\ 2 u\ +\ 2\ =\ 0\,, \qquad
  u(-1)\ =\ u(1)\ =\ 0\,.
\end{equation}
This BVP has the exact solution
\begin{equation}\label{eq:exa}
  u(x)\ =\ 1\ -\ \frac{\sinh(2)}{\sinh(3)}\;\ue^{x}\ -\ \frac{\sinh(1)}{\sinh(3)}\;\ue^{-2 x}\,, \qquad x\ \in\ [-1,\,1]\,.
\end{equation}
We shall seek for the numerical solution to \eqref{eq:bvp} in the form of a truncated \textsc{Tchebyshev} expansion:
\begin{equation*}
  u(x)\ \approx\ \sum_{k\, =\, 0}^{4} a_k\, T_k(x)\,, \qquad \forall k:\ a_k\ \in\ \R
\end{equation*}
Spectral coefficients $\{a_k\}_{k = 0}^{4}$ have to be determined. The last expansion is substituted into the governing equation \eqref{eq:bvp}. There is no reason that \eqref{eq:bvp} will be satisfied identically in every point of $\I$. Hence, we can measure the residual
\begin{equation*}
  \Rr(x)\ =\ \bigl(u_{xx}\ +\ u_x\ -\ 2 u\ +\ 2\bigr)(x)\ \rightsquigarrow\ 0\,.
\end{equation*}
The enforcing of boundary conditions $u(-1)\ =\ u(1)\ =\ 0$ leads to two additional relations on spectral coefficients:
\begin{align*}
  a_0\ +\ a_1\ +\ a_2\ +\ a_3\ +\ a_4\ &=\ 0\,, \\
  a_0\ -\ a_1\ +\ a_2\ -\ a_3\ +\ a_4\ &=\ 0\,.
\end{align*}
So, we have two relations coming from boundary conditions and we have five degrees of freedom $\{a_k\}_{k = 0}^{4}$. It means that we have to impose three additional conditions to determine uniquely all spectral coefficients. Different approaches prescribe different numerical recipes.

\subsubsection{Tau--Lanczos}

We require that the residual $\Rr(x)$ be orthogonal to the first three basis functions $T_{0,1,2}(x)$. It gives us three additional relations:
\begin{equation*}
  \langle\Rr,\, T_k\rangle\ \equiv\ \int_{-1}^{1} \frac{\Rr(x)\, T_k(x)}{\sqrt{1\ -\ x^2}}\;\ud x\ =\ 0\,, \qquad k\ =\ 0,\, 1,\, 2\,.
\end{equation*}

\subsubsection{Galerkin}

We recombine the \textsc{Tchebyshev} polynomials to form a different basis which satisfies identically the boundary conditions:
\begin{align*}
  \phi_0(x)\ &\eqdef\ \bigl(T_2\ -\ T_0\bigr)(x)\,, \\
  \phi_1(x)\ &\eqdef\ \bigl(T_3\ -\ T_1\bigr)(x)\,, \\
  \phi_2(x)\ &\eqdef\ \bigl(T_4\ -\ T_0\bigr)(x)\,. \\
\end{align*}
Then we require that the residual $\Rr(x)$ is orthogonal to the new basis functions:
\begin{equation*}
  \langle\Rr,\, \phi_k\rangle\ \equiv\ \int_{-1}^{1} \frac{\Rr(x)\, \phi_k(x)}{\sqrt{1\ -\ x^2}}\;\ud x\ =\ 0\,, \qquad k\ =\ 0,\, 1,\, 2\,.
\end{equation*}

\subsubsection{Collocation}

This approach is particularly simple. We require that $\Rr(x_k)\ \equiv\ 0$ in three interior \textsc{Tchebyshev} points $x_k\ =\ \cos\bigl(\frac{\pi k}{4}\bigr)$, $k = 1,\, 2,\, 3$.

\begin{figure}
  \centering
  \subfigure[Exact and numerical solutions]{
  \includegraphics[width=0.485\textwidth]{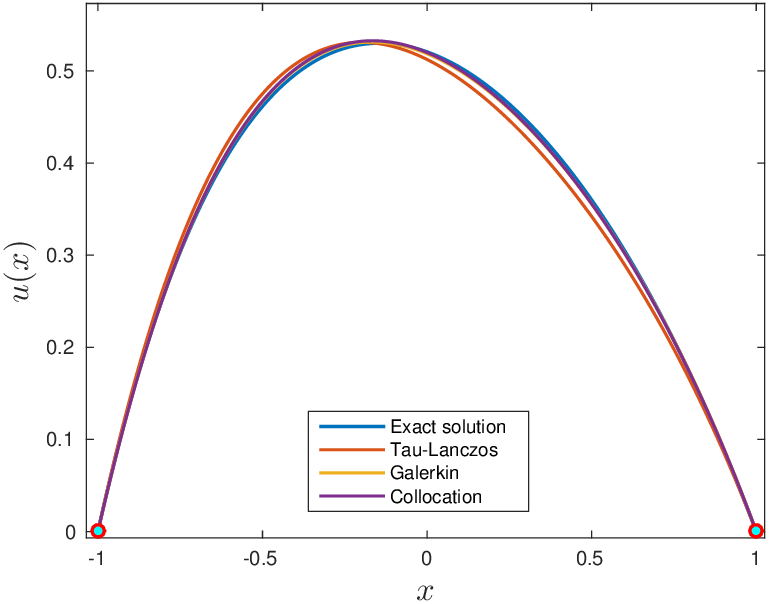}}
  \subfigure[Zoom on a portion of domain]{
  \includegraphics[width=0.485\textwidth]{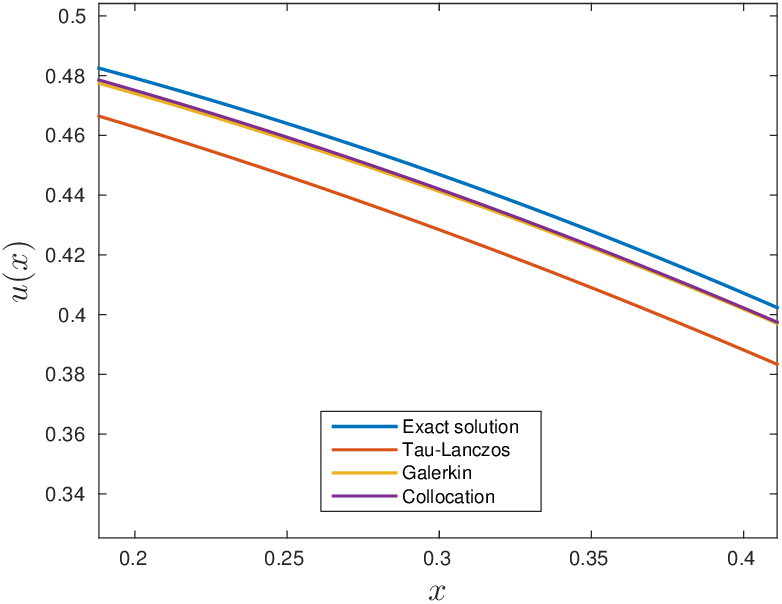}}
  \caption{\small\em Comparison of various numerical approaches to determine the spectral expansion coefficients for a linear BVP \eqref{eq:bvp}.}
  \label{fig:sols}
\end{figure}

\bigskip
\paragraph*{Conclusions.} The comparison of three numerical solutions with the exact one \eqref{eq:exa} is shown in Figure~\ref{fig:sols}. In this Figure~\ref{fig:sols} one can see that the collocation\footnote{The advantage of the collocation approach for nonlinear (and variable coefficients) problems becomes even more flagrant.} and \textsc{Galerkin} methods provide nearly identical\footnote{With a little advantage towards the collocation method. However, this conclusion is not general at all. It is based only on this particular problem.} numerical solutions which are closer to the exact solution than the prediction of the Tau method. However, in general the Author is impressed by the performance of spectral methods --- with only five degrees of freedom $\{a_k\}_{k = 0}^{4}$ we capture very well the global behaviour of the solution \eqref{eq:exa} in every point of the interval $\I$. Perhaps, the last point has to be emphasized again: the (pseudo-)spectral methods provide the numerical value of the solution in the \emph{whole domain}, not only in collocation or grid points.

\section{Application to heat conduction}
\label{sec:heat}

In this Section we show how to apply the \textsc{Fourier}-type spectral methods to simulate some simple and not so simple diffusion processes.

\subsection{An elementary example}

\begin{figure}
  \centering
  \includegraphics[width=0.75\textwidth]{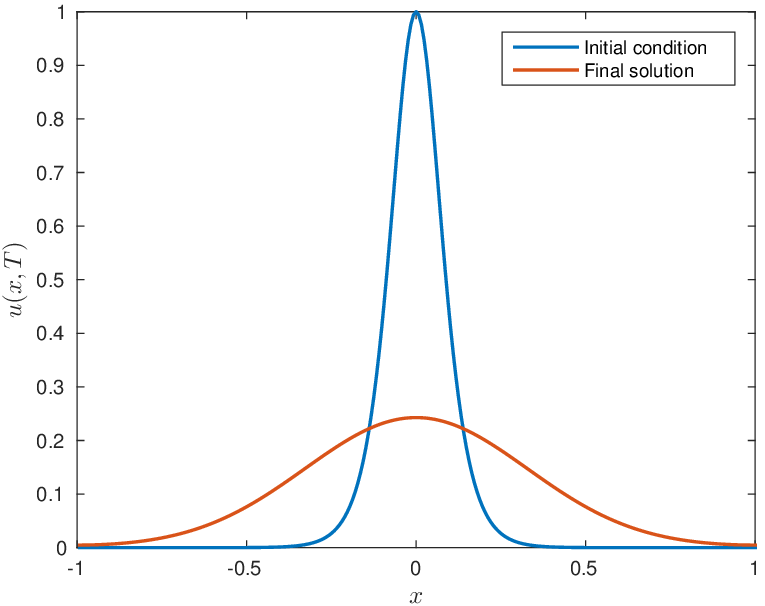}
  \caption{\small\em Spectral solution to the linear heat equation \eqref{eq:heatF} at $t\ =\ 5.0\,$, with $\nu\ =\ 10^{-2}$ and the initial condition is $u_0(x)\ =\ \sech^2(10\,x)\,$.}
  \label{fig:heatF}
\end{figure}

As the simplest example, consider the linear heat equation
\begin{equation}\label{eq:heatF}
  u_t\ =\ \nu\,u_{xx}\,, \qquad x\ \in\ \R\,,
\end{equation}
completed with the initial condition
\begin{equation*}
  u(x,\,0)\ =\ u_0(x)\,.
\end{equation*}
We consider this equation on the whole line $\R$ for the sake of simplicity. However, \emph{in silico} it will be truncated to a periodic interval (in the example below it will be $\I\ =\ [-1,\,1]\,$). Let us apply the \textsc{Fourier} integral transform to the both sides of equation \ref{eq:heatF}:
\begin{equation}\label{eq:heatT}
  \od{\hat{u}(k,\,t)}{t}\ =\ -\nu\,k^2\,\hat{u}(k,\,t)\,,
\end{equation}
where the forward $\hat{u}(k,\,t)\ =\ \F\{u(x,\,t)\}$ and inverse $u(x,\,t)\ =\ \F^{-1}\{\hat{u}(k,\,t)\}$ integral \textsc{Fourier} transforms are defined below in \eqref{eq:forw} and \eqref{eq:inv} correspondingly. The transformed heat equation \eqref{eq:heatT} can be regarded as a linear ODE. Its solution can be readily obtained:
\begin{equation*}
  \hat{u}(k,\,t)\ =\ \hat{u}_0(k)\,\ue^{-\nu\,k^2\,t}\,, \qquad
  \hat{u}_0(k)\ =\ \F\{u_0(x)\}\,.
\end{equation*}
Consequently, we possess an analytical solution to the heat equation \eqref{eq:heatF} in \textsc{Fourier} space. In order to obtain the solution in physical space, an inverse \textsc{Fourier} transform has to be computed. The exponential decay of \textsc{Fourier} coefficients $\hat{u}(k,\,t)$ (except the zero mode $k\ =\ 0$) ensures that the solution $u(x,\,t)$ becomes infinitely smooth $C^\infty$ for $t\ >\ 0\,$. For example, the initial condition and corresponding solution at $t\ =\ T\ =\ 5\ \s$ are represented in Figure~\ref{fig:heatF}. The \textsc{Matlab} code used to generate Figure~\ref{fig:heatF} is provided below as well.

\begin{lstlisting}
l  = 1.0;    % half-length of the domain
N  = 256;    % number of Fourier modes
dx = 2*l/N;  % distance between two collocation points
x  = (1-N/2:N/2)*dx; % physical space discretization
nu = 0.01;   % diffusion parameter
T  = 5.0;    % time where we compute the solution

dk = pi/l;                  % discretization step in Fourier space
k  = [0:N/2 1-N/2:-1]*dk;   % vector of wavenumbers
k2 = k.^2;                  % almost 2nd derivative in Fourier space

u0     = sech(10.0*x).^2;   % initial condition
u0_hat = fft(u0);           % Its Fourier transform

% and the solution at final time:
uT     = real(ifft(exp(-nu*k2*T).*u0_hat));
\end{lstlisting}

\subsection{A less elementary example}

The numerical code above is based on the knowledge of the analytical solution to ODE \eqref{eq:heatT} in the \textsc{Fourier} space. It is unnecessary to say that an analytical solution is available in very simple situations only. Consequently, we assume that the resulting ODE system in the physical (\textsc{Fourier}) space is formally solved (\ie~advanced in time) by the semigroup operator $\S(t)$ ($\Sh(t)\ \equiv\ \F\cdot\S(t)\cdot\F^{\,-1}$ correspondingly):
\begin{equation*}
  u(x,\,t)\ \equiv\ \S(t)\cdot\bigl[\,u_0(x)\,\bigr]\,.
\end{equation*}
In practice this semigroup operator is realized using numerical time-marching techniques described briefly in Appendix~\ref{app:odes}. So, a more general \textsc{Fourier} spectral algorithm is
\begin{enumerate}
  \item Decompose the initial condition:
  \begin{equation*}
    \hat{u}_0(k)\ =\ \F\bigl\{u_0(x)\bigr\}\,,
  \end{equation*}
  \item Advance in time:
  \begin{equation*}
    \hat{u}_k(t)\ =\ \Sh(t)\bigl[\,\hat{u}_0(k)\,\bigr]\,,
  \end{equation*}
  \item Synthesize:
  \begin{equation*}
    u(x,\,t)\ =\ \F^{\,-1}\bigl\{\hat{u}_k(t)\bigr\}\,.
  \end{equation*}
\end{enumerate}
This algorithm is based on the fact that the following diagram commutes:
\begin{equation*}
\begin{tikzcd}[row sep = large, column sep = large]
  u(x,\,0) \arrow{d}[swap]{\F} \arrow{r}[blue]{\S(t)}
  & u(x,\,t) \\
  \hat{u}_0(k) \arrow{r}[blue]{\Sh(t)}
  & \hat{u}_k(t) \arrow{u}{\F^{\,-1}}
\end{tikzcd}  
\end{equation*}

\subsection{A real-life example}

In this Section we shall consider a realistic model (to be honest we take a slightly simplified version), which was proposed in \cite{Mendes2005} to predict heat and moisture transfer through the walls. This model is used in real-world Civil Engineering applications such as the code \textsc{Domus}. So, the model we consider in this Section reads (for simplicity we consider the 1D case):
\begin{align}\label{eq:nat1}
  \pd{\theta}{t}\ &=\ \pd{}{x}\Bigl(\D_\theta\,\pd{\theta}{x}\ +\ \D_T\,\pd{T}{x}\Bigr)\,, \\
  \rho c_m\,\pd{T}{t}\ &=\ -\pd{q}{x}\ -\ L(T)\cdot\pd{j_v}{x}\,,\label{eq:nat2}
\end{align}
where $\theta$ is moisture volumetric content (\ie~moisture density) and $T$ is the temperature. The mass transport coefficients $\D_\theta(\theta,\, T)$ and $\D_T(\theta,\, T)$ may depend nonlinearly on the solution $\bigl(\theta(x,\,t),\, T(x,\,t)\bigr)$. Here we assume for simplicity that the mass density $\rho$ and the specific moisture heat $c_m$ are some positive constants\footnote{In more realistic modelling $c_m(\theta)$ depends also on the moisture density $\theta$. It does not pose any problems to take it into account in the numerical scheme described in Section~\ref{sec:heat}.}. Finally, the heat flux $q$ and vapor flow $j_v$ can be expressed as
\begin{align*}
  q\ &=\ -\lambda(\theta,\,T)\,\pd{T}{x}\,, \\
  j_v\ &=\ -\V_\theta(\theta,\,T)\,\pd{\theta}{x}\ -\ \V_T(\theta,\,T)\,\pd{T}{x}\,.
\end{align*}

In \textsc{Fourier}-type pseudo-spectral methods we usually work with \textsc{Fourier} coefficients. Consequently, we apply the \textsc{Fourier} transform to both sides of equations \eqref{eq:nat1}, \eqref{eq:nat2}:
\begin{align}\label{eq:nat3}
  \od{\thetah}{t}\ &=\ \ui k\,\F\Bigl\{\D_\theta\,\pd{\theta}{x}\ +\ \D_T\,\pd{T}{x}\Bigr\}\,, \\
  \rho c_m\,\od{\Th}{t}\ &=\ \ui k\,\F\Bigl\{\lambda(\theta,\,T)\,\pd{T}{x}\Bigr\}\ -\ \F\Bigl\{L(T)\cdot\pd{j_v}{x}\Bigr\}\,,\label{eq:nat4}
\end{align}
where we introduced some notations for the \textsc{Fourier} transform $\F(\cdot)$:
\begin{equation}\label{eq:forw}
  \thetah(k,\,t)\ \eqdef\ \F\Bigl\{\theta(x,\,t)\Bigr\}\ =\ \int_{-\infty}^{+\infty} \theta(x,\,t)\,\ue^{\ui k x}\;\ud x\,,
\end{equation}
where $k \in \R$ is the wave number. Inversely we have
\begin{equation}\label{eq:inv}
  \theta(x,\,t)\ =\ \F^{\,-1}\Bigl\{\thetah(k,\,t)\Bigr\}\ =\ \frac{1}{2\pi}\,\int_{-\infty}^{+\infty} \thetah(k,\,t)\,\ue^{-\ui k x}\;\ud k\,.
\end{equation}
Now, system \eqref{eq:nat3}, \eqref{eq:nat4} can be considered as a coupled nonlinear system of Ordinary Differential Equations (ODEs) (not PDEs!) for the \textsc{Fourier} coefficients of solutions $\bigl(\theta(x,\,t),\, T(x,\,t)\bigr)$. The numerical methods for systems of ODEs will be explained in a separate course. As general fundamental references on this topic we can recommend \cite{Hairer2002, Hairer2009}. The Author even suggests to employ a well-documented ready-to-use ODE library such as \cite{Shampine1997}, for example. Now we have to explain how to evaluate the right hand side of equations \eqref{eq:nat3}, \eqref{eq:nat4} when only the \textsc{Fourier} coefficients are available. The recipe is very simple:
\begin{itemize}
  \item All linear operations (\eg additions, subtractions, multiplication by a scalar) can be equally made in \textsc{Fourier} or in a real spaces. The choice has to be done in order to minimize the number of FFT operations
  \begin{equation*}
    \F\bigl\{\alpha\theta\ +\ \beta T\bigr\}\ \equiv\ \alpha\thetah\ +\ \beta\Th
  \end{equation*}
  \item All nonlinear products are made in the real space and then we transfer the result back to the \textsc{Fourier} space using one FFT, \eg
  \begin{equation*}
    \F\bigl\{\theta(x,\,t)\cdot T(x,\,t)\bigr\}\ =\ \F\Bigl\{\F^{\,-1}\bigl\{\thetah(k,\,t)\bigr\}\cdot\F^{\,-1}\bigl\{\Th(k,\,t)\bigr\}\Bigr\}
  \end{equation*}
  \item All spatial derivatives are computed only in \textsc{Fourier} space, as
  \begin{equation*}
    \pd{^n}{x^n}\bigl[\,\cdot\,\bigr]\quad \leftrightarrows\quad \F^{\,-1}\,\Bigl\{(\ui k)^n\,\bigl[\,\cdot\,\bigr]\Bigr\}
  \end{equation*}
\end{itemize}
The Reader can notice that above we used \textsc{Fourier} integrals. It comes from the mathematical tradition. In computer implementations one has to take a finite interval, periodize it and use direct and inverse Discrete \textsc{Fourier} Transforms (DFTs) instead of $\F\{\cdot\}$ and $\F^{\,-1}\{\cdot\}$ respectively.

\begin{figure}
  \centering
  \includegraphics[width=0.99\textwidth]{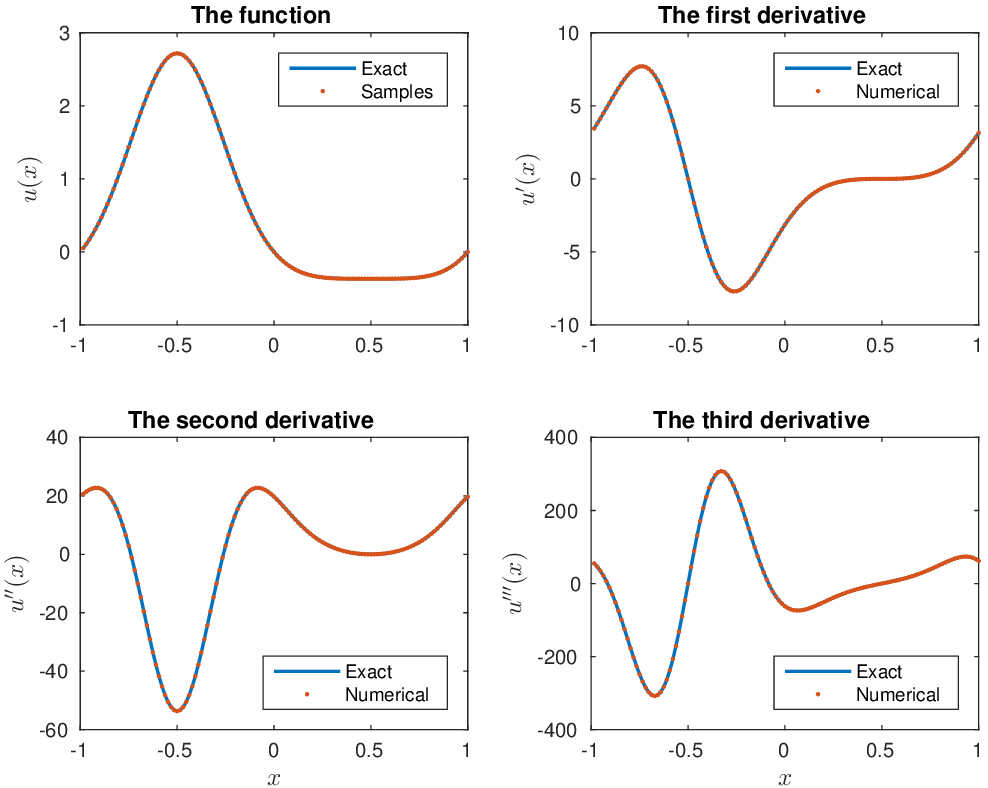}
  \caption{\small\em Numerical differentiation of a periodic function \eqref{eq:func} using the \textsc{Fourier} collocation spectral method. We use $128$ collocation points. The agreement is perfect up to graphic accuracy.}
  \label{fig:diffF}
\end{figure}

The computation of derivatives using the \textsc{Fourier} collocation spectral method is illustrated below on the following periodic function (with period $2$):
\begin{equation}\label{eq:func}
  u(x)\ =\ \sin\bigl(\pi(x\ +\ 1)\bigr)\,\ue^{\sin\left(\pi(x\ +\ 1)\right)}\,, \qquad x\ \in\ [-1,\,1]\,.
\end{equation}
The first three derivatives of this function can be readily computed using any symbolic computation software (the Author used \textsc{Maple}):
\begin{equation*}
  u^{\,\prime}(x)\ =\ \pi\,\cos\bigl(\pi(x\ +\ 1)\bigr)\Bigl(1\ +\ \sin\bigl(\pi(x\ +\ 1)\bigr)\Bigr)\,\ue^{\sin\left(\pi(x\ +\ 1)\right)}\,,
\end{equation*}
\begin{equation*}
  u^{\,\prime\prime}(x)\ =\ \pi^2\,\Bigl\{\cos^2\bigl(\pi(x\ +\ 1)\bigr)\,\Bigl(\sin\bigl(\pi(x\ +\ 1)\bigr)\ +\ 3\Bigr)\ -\ \sin\bigl(\pi(x\ +\ 1)\bigr)\ -\ 1\Bigr\}\,\ue^{\sin\left(\pi(x\ +\ 1)\right)}\,,
\end{equation*}
\begin{equation*}
  u^{\,\prime\prime\prime}(x)\ =\ \pi^3\,\cos\bigl(\pi(x\ +\ 1)\bigr)\,\Bigl\{\cos^2\bigl(\pi(x\ +\ 1)\bigr)\,\Bigl(\sin\bigl(\pi(x\ +\ 1)\bigr)\ +\ 6\Bigr)\ -\ 7\,\sin\bigl(\pi(x\ +\ 1)\bigr)\ -\ 4\Bigr\}\,.
\end{equation*}
The expressions above are used to assess the accuracy of the computed approximations in collocation points. The analytical derivatives with computed ones (red dots) are shown in Figure~\ref{fig:diffF} for $N\ =\ 128$ collocation points. To the graphical accuracy the results are indistinguishable. That is why we computed also the discrete $\ell_\infty$ norms to compute the relative errors:
\begin{equation*}
  \eps_N^{(p)}\ =\ \frac{\norm{u_{\mathrm{num}}\ -\ u_{\mathrm{exact}}}_{\infty}}{\norm{u_{\mathrm{exact}}(x_i)}_{\infty}}\,.
\end{equation*}
For $N\ =\ 32$ \textsc{Fourier} modes we obtain the following numerical results:
\begin{align*}
  \eps_{32}^{(1)}\ &\approx\ 1.5\times 10^{-15}\,, \\
  \eps_{32}^{(2)}\ &\approx\ 8.4\times 10^{-15}\,, \\
  \eps_{32}^{(3)}\ &\approx\ 4.7\times 10^{-14}\,.
\end{align*}
The \textsc{Matlab} code used to generate these results and Figure~\ref{fig:diffF} is provided below.
\begin{lstlisting}
l  = 1.0;    % half-length of the domain
N  = 128;    % number of Fourier modes
dx = 2*l/N;  % distance between two collocation points
x  = (1-N/2:N/2)*dx; % physical space discretization

dk = pi/l;
k  = [0:N/2 1-N/2:-1]*dk; % vector of wavenumbers

arg  = pi*(x + 1);
sar  = sin(arg);
car  = cos(arg); car2 = car.*car;
esa  = exp(sar);
u    = sar.*esa;
uhat = fft(u);

% numerical derivatives:
up   = real(ifft(1i*k.*uhat));
upp  = real(ifft(-k.^2.*uhat));
uppp = real(ifft(-1i*k.^3.*uhat));

% exact derivatives:
pi2  = pi*pi; pi3 = pi2*pi;
u1   = pi*car.*(1 + sar).*esa;
u2   = pi2*(car2.*(sar + 3) - 1 - sar).*esa;
u3   = pi3*car.*(car2.*(sar + 6) - 4 - 7*sar).*esa;
\end{lstlisting}
Above we exploit the property that the \textsc{Fourier} transform is a `change a variables' where the differentiation operator $\partial_x$ becomes diagonal.

\section{Indications for further reading}
\label{sec:indic}

If you got interested in the beautiful topic of pseudo-spectral methods, you can find more information in the following books:

\begin{itemize}

  \item The following book by N.~\textsc{Trefethen}\footnote{Lloyd Nick~\textsc{Trefethen} (1955 -- 20..), an American/British numerical analyst.} has two major advantages: ({\it i}) conciseness and ({\it ii}) collection of \textsc{Matlab} programs which comes along. I would recommend it as the first reading on pseudo-spectral methods. At least you will learn how to program efficiently in \textsc{Matlab} (or in \textsc{Octave}, \textsc{Scilab}, \etc) \cite{Trefethen2000}:
  \begin{itemize}
    \item \textsc{Trefethen}, L. N. (2000). \textit{Spectral methods in MatLab}. Society for Industrial and Applied Mathematics, Philadelphia, PA, USA.
  \end{itemize}
  \smallskip

  \item The book by J.~\textsc{Boyd}\footnote{John~\textsc{Boyd} (19.. -- 20..), an American numerical analyst, meteorologist and occasional science fiction writer.} is probably the most exhaustive one. It covers many topics and applications of spectral methods. The material is presented as a collection of tricks. For example, it is one of seldom books where (semi-)infinite domains and spherical geometries are covered. I am not sure that after reading this book you will know how to program the pseudo-spectral methods, but you will have a broad of view of possible issues and how to address them \cite{Boyd2000}:
  \begin{itemize}
    \item \textsc{Boyd}, J. P. (2000). \textit{Chebyshev and Fourier Spectral Methods}. (Dover Publications, New York) (2\up{nd} Ed.).
  \end{itemize}
  \smallskip

  \item This book is probably my favourite one. It represents a good balance between the theory and practice and this book arose as an extended version of a previously published review paper. The present lecture notes were inspired in part on this book as well \cite{Fornberg1996}:
  \begin{itemize}
    \item \textsc{Fornberg}\footnote{Bengt~\textsc{Fornberg} (19.. -- 20..), a Swedish/American numerical analyst.}, B. (1996). \textit{A practical guide to pseudospectral methods}. Cambridge: Cambridge University Press.
  \end{itemize}

  \item The Author discovered this book only recently. So, he has still to read it, but from the first sight I can already recommend it:
  \begin{itemize}
    \item \textsc{Peyret}, R. (2002). \textit{Spectral Methods for Incompressible Viscous Flow}. Springer-Verlag New York Inc.
  \end{itemize}

  \item Finally, I discovered also these very clear and instructive Lecture notes \cite{Uecker2009}. Lectures were delivered in 2009 at the International Summer School `\textit{Modern Computational Science}' in Oldenburg, Germany. They can be freely downloaded at the following URL address:
  \begin{itemize}
    \item \url{http://www.staff.uni-oldenburg.de/hannes.uecker/pre/030-mcs-hu.pdf}
  \end{itemize}

\end{itemize}

In general, the Author's suggestions (according to his personal taste and vision) for scientific literature and software are collected in a single document, which is continuously expanded and completed:

\begin{center}
  \url{https://github.com/dutykh/libs/}
\end{center}
\smallskip
\begin{center}
  \textsc{Good luck and have a nice reading!}
\end{center}
\begin{center}
  \Large
  $\star$ $\star$ $\star$
\end{center}


\appendix


\section{Some identities involving Tchebyshev polynomials}
\label{app:tcheb}

\begin{figure}
  \centering
  \includegraphics[width=0.99\textwidth]{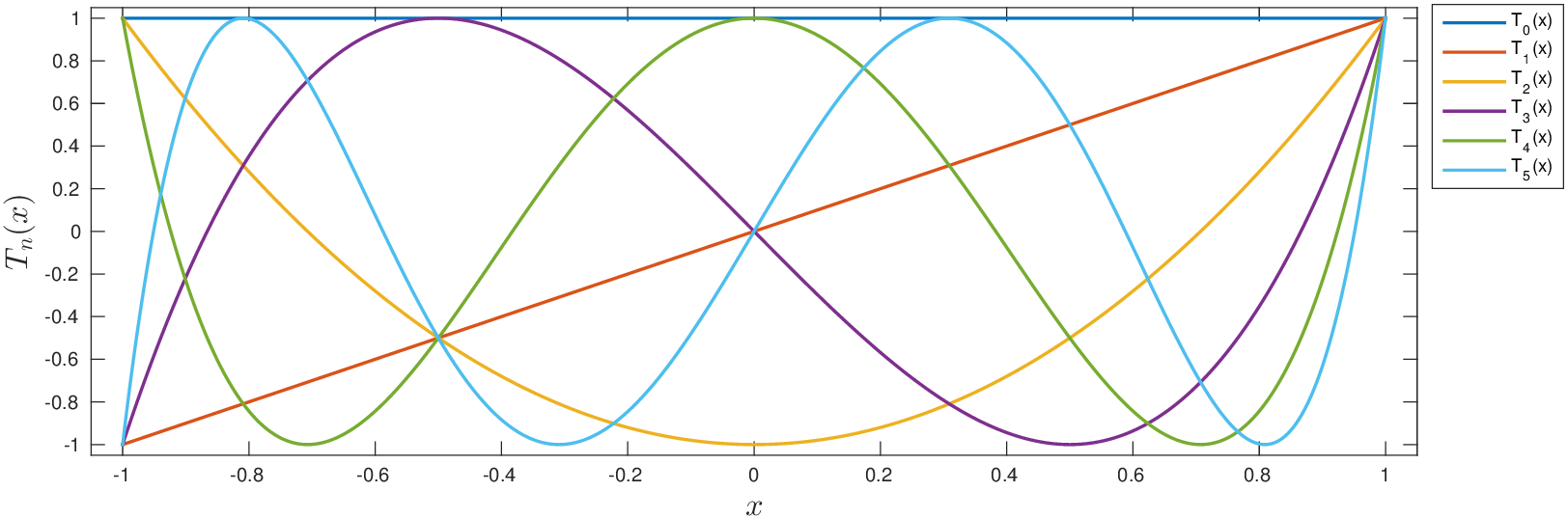}
  \caption{\small\em The first five \textsc{Tchebyshev} polynomials. Notice that $T_n(-1)\ =\ (-1)^n$ and $T_n(1)\ \equiv\ 1$.}
  \label{fig:tcheb}
\end{figure}

\textsc{Tchebyshev} polynomials $\{T_n(x)\}_{n = 0}^{\infty}$ form a weighted orthogonal system on the segment $[-1,\, 1]$:
\begin{equation*}
  \langle\, T_n,\, T_m\,\rangle\ \equiv\ \int_{-1}^{1}\,\frac{T_n(x) \cdot T_m(x)}{\sqrt{1\ -\ x^2}}\;\ud x\ =\ \begin{dcases}
    0,\quad m\ \neq\ n\,, \\
    \pi,\quad m\ =\ n\ =\ 0\,, \\
    \frac{\pi}{2},\quad m\ =\ n\ >\ 0\,.
  \end{dcases}
\end{equation*}
The first few \textsc{Tchebyshev} polynomials are
\begin{align*}
  T_0(x)\ &=\ 1\,, \\
  T_1(x)\ &=\ x\,, \\
  T_2(x)\ &=\ 2\, x^2\ -\ 1\,, \\
  T_3(x)\ &=\ 4\, x^3\ -\ 3\, x\,, \\
  T_4(x)\ &=\ 8\, x^4\ -\ 8\, x^2\ +\ 1\,, \\
  T_5(x)\ &=\ 16\, x^5\ -\ 20\, x^3\ +\ 5\, x\,.
\end{align*}
They are represented graphically in Figure~\ref{fig:tcheb}. Higher order \textsc{Tchebyshev} polynomials can be constructed using the three-term recursion relation:
\begin{equation*}
  T_{n+1}(x)\ =\ 2\, x\, T_n(x)\ -\ T_{n-1}(x)\,.
\end{equation*}
The following \textsc{Matlab} code realizes this idea in practice:
\begin{lstlisting}
function P = Chebyshev(n)
    P = cell(n,1);
    if (n == 1)
        P{1} = 1;
        return;
    end % if ()

    P{1} = 1;
    P{2} = [1; 0];
    for j=3:n
        P{j}      = [2*P{j-1}; 0];
        P{j}(3:j) = P{j}(3:j) - P{j-2};
    end % for j
    
end % Chebyshev()
\end{lstlisting}
\bigskip

Then, $n$\up{th} \textsc{Tchebyshev} polynomial can be evaluated using the standard \textsc{Matlab}'s function \texttt{polyval()}, \eg
\begin{lstlisting}
n  = 5;
P  = Chebyshev(n+1);
x  = linspace(-1, 1, 1000);
Tn = polyval(P{n+1}, x);
\end{lstlisting}
There is an explicit expression for the $n$\up{th} polynomial:
\begin{equation}\label{eq:cos}
  T_n(x)\ =\ \cos\bigl(\,n\theta(x)\,\bigr)\,, \qquad \theta(x)\ \eqdef\ \arccos x\,.
\end{equation}
The first derivatives of \textsc{Tchebyshev} polynomials can be constructed recursively as well:
\begin{equation}\label{eq:der}
  \frac{T_{n+1}^{\,\prime}(x)}{n\ +\ 1}\ =\ \frac{T_{n-1}^{\,\prime}(x)}{n\ -\ 1}\ +\ 2\, T_n(x)\,.
\end{equation}
Otherwise, the first derivative $T_n^{\,\prime}(x)$ can be found from the following relation:
\begin{equation*}
  (1\ -\ x^2)\,T_n^{\,\prime}(x)\ =\ -n\, x\, T_n(x)\ +\ n\, T_{n - 1}(x)\,.
\end{equation*}
\textsc{Tchebyshev} polynomial $T_n(x)$ satisfies the following linear second order differential equation with non-constant coefficients:
\begin{equation*}
  (1\ -\ x^2)\, T_n^{\,\prime\prime}\ -\ x\, T_n^{\,\prime}\ +\ n^2\, T_n\ =\ 0\,.
\end{equation*}
Finally, zeros of the $n$\up{th} \textsc{Tchebyshev} polynomial $T_n(x)$ are located at
\begin{equation*}
  x_k^{0}\ =\ \cos\,\Bigl(\,\frac{2k - 1}{2n}\,\pi\,\Bigr)\,, \qquad k\ =\ 1,\, 2,\,\ldots,\, n,
\end{equation*}
and extrema at
\begin{equation*}
  x_k^{\mathrm{ext}}\ =\ \cos\,\Bigl(\,\frac{\pi k}{n}\,\Bigr)\,, \qquad k\ =\ 0,\, 1,\,\ldots,\, n\,.
\end{equation*}

When implementing \textsc{Tchebyshev}-type pseudo-spectral methods for nonlinear problems, one can use the expression of the product of two \textsc{Tchebyshev} polynomials in terms of higher and lower degree polynomials:
\begin{equation*}
  T_m(x)\,T_n(x)\ =\ \frac{1}{2}\;\bigl(\,T_{n + m}(x)\ +\ T_{n - m}(x)\,\bigr)\,,\qquad \forall n\ \geq\ m\ \geq\ 0\,.
\end{equation*}
To finish this Appendix we give the \emph{generating function} for \textsc{Tchebyshev} polynomials:
\begin{equation*}
  \frac{1\ -\ x\,z}{1\ -\ 2\, x\, z\ +\ z^2}\ =\ \sum_{n\, =\, 0}^{\infty}\, T_n(x)\, z^n\,.
\end{equation*}

In order to satisfy the requirement (2) from Section~\ref{sec:basis}, we provide here the relations which allow to re-express the derivatives of the \textsc{Tchebyshev} expansion in terms of \textsc{Tchebyshev} polynomials again. Namely, consider a truncated series expansion of a function $u(x)$ in \textsc{Tchebyshev} polynomials:
\begin{equation*}
  u(x)\ =\ \sum_{k\,=\,0}^{N} v_k\,T_k(x)\,.
\end{equation*}
Imagine that we want to compute the first derivative of the expansion above (it is always possible, since the sum is finite and we can differentiate term by term):
\begin{equation*}
  u^{\,\prime}(x)\ =\ \sum_{k\,=\,0}^{N} v_k\,\od{T_k(x)}{x}\,.
\end{equation*}
Now we re-expand the derivative $u^{\,\prime}(x)$ in the same basis functions:
\begin{equation*}
  u^{\,\prime}(x)\ =\ \sum_{k\,=\,0}^{N} v_k^{\,\prime}\,T_k(x)\,.
\end{equation*}
The main question is how to re-express the coefficients $\{v_k^{\,\prime}\}_{k=0}^{N}$ in terms of coefficients $\{v_k\}_{k=0}^{N}\,$? This goal is achieved using basically the recurrence relation \eqref{eq:der} (even if it does not jump into the eyes). Thanks to it we have an explicit relation between the coefficients:
\begin{equation*}
  v_k^{\,\prime}\ =\ \frac{2}{\delta_k}\;\sum_{\stackrel{j\, =\, k+1}{j+k\ \mathrm{odd}}}^{N}\, j\,v_j\,, \qquad j\ =\ 0,\,1,\,\ldots,N-1\,,
\end{equation*}
where
\begin{equation*}
  \delta_k\ =\ \begin{dcases}
    \ 2,& \ k\,=\,0\,, \\
    \ 0,& \ k\,\neq\,0\,.
  \end{dcases}
\end{equation*}
A similar `trick' can be made for the 2\up{nd} derivative as well:
\begin{equation*}
  u^{\,\prime\prime}(x)\ =\ \sum_{k=0}^{N} v_k\,\od{^2T_k(x)}{x^2}\ =\ \sum_{k=0}^{N} v_k^{\,\prime\prime}\,T_k(x)\,.
\end{equation*}
The connection between coefficients $\{v_k^{\,\prime\prime}\}_{k=0}^{N}$ and $\{v_k\}_{k=0}^{N}$ is given by the following explicit formula:
\begin{equation*}
  v_k^{\,\prime\prime}\ =\ \frac{1}{\delta_k}\;\sum_{\stackrel{j = k+2}{j+k\ \mathrm{even}}}^{N}\, j\,(j^2\ -\ k^2)\,v_j\,, \qquad j\ =\ 0,\,1,\,\ldots,N-2\,.
\end{equation*}
Finally, in order to construct the spectral coefficients for the $n$\up{th} derivative, one can use the following recurrence relation (also stemming from \eqref{eq:der}):
\begin{equation*}
  \delta_{k-1}\,v_{k-1}^{(n)}\ =\ v_{k+1}^{(n)}\ +\ 2\,k\,v_k^{(n-1)}\,, \qquad k\ \geq\ 1\,,
\end{equation*}
which has to be completed by starting value $v_N^{\,\prime}\ \equiv\ 0\,$.

\begin{remark}
Similar identities exist for the family of \textsc{Jacobi} polynomials $\Pl_n^{\,\alpha,\,\beta}(x)$ as well. However, they are more complicated due to the presence of two arbitrary parameters $\alpha,\, \beta\ >\ -1$. The family of \textsc{Tchebyshev} polynomials is a particular case of \textsc{Jacobi} polynomials when $\alpha\ =\ \beta\ =\ -\frac{1}{2}$.
\end{remark}

\subsection{Compositions of Tchebyshev polynomials}

We would like to show here another interesting formula involving \textsc{Tchebyshev} polynomials:
\begin{theorem}[Composition formula]
If $T_n(x)$ and $T_m(x)$ are Tchebyshev polynomials ($m, n\ \geq\ 0$), then
\begin{equation}\label{eq:comp}
  \bigl(T_m\circ T_n\bigr)\,(x)\ \equiv\ T_m\bigl(T_n\,(x)\bigr)\ =\ T_{mn}(x)\,.
\end{equation}
\end{theorem}

\begin{proof}
An \emph{eleven} (!) pages combinatorial proof of this result can be found in \cite{Benjamin2010}. Here we shall prove the composition formula \eqref{eq:comp} in one line by using some rudiments of complex variables. Let us introduce the variable $z\ \eqdef\ \ue^{\ui\theta}\ \in\ \circlearrowleft$\footnote{Symbol $\circlearrowleft$ denotes the unit circle $\S_1(\vO)$ on the complex plane $\C$.} such that
\begin{equation*}
  x\ =\ \frac{1}{2}\;\Bigl(z\ +\ \frac{1}{z}\Bigr)\, \quad \leftrightarrows \quad x\ =\ \cos\theta\,.
\end{equation*}
Then, by \eqref{eq:cos} we have a complex representation of the $n$\up{th} \textsc{Tchebyshev} polynomial:
\begin{equation*}
  T_n\,(x)\ \equiv\ T_n\Bigl[\half\,\bigl(z\ +\ z^{\,-1}\bigr)\Bigr]\ =\ \frac{1}{2}\;\Bigl(z^{\,n}\ +\ \frac{1}{z^{\,n}}\Bigr)\,.
\end{equation*}
Now, we have all the elements to show the main result:
\begin{multline*}
T_m\biggl[\,T_n\Bigl(\half\,\bigl(z\ +\ z^{\,-1}\bigr)\Bigr)\,\biggr]\ =\ T_m\Bigl[\,\half\,\bigl(z^{\,n}\ +\ z^{\,-n}\bigr)\Bigr]\ =\ \frac{1}{2}\;\Bigl((z^{\,n})^{\,m}\ +\ \frac{1}{(z^{\,n})^{\,m}}\Bigr)\ =\\
\frac{1}{2}\;\Bigl(z^{\,mn}\ +\ \frac{1}{z^{\,mn}}\Bigr) =\ T_{mn}\Bigl[\half\,\bigl(z\ +\ z^{\,-1}\bigr)\Bigr]\ \equiv\ T_{mn}\,(x)\,.
\end{multline*}
\end{proof}


\section{Trefftz method}
\label{app:trefftz}

In this Appendix I would like to give the flavour of the so-called \textsc{Trefftz}\footnote{Erich Emmanuel \textsc{Trefftz} (1888 -- 1937), a German Mathematician and Mechanical Engineer.} methods, which remain essentially unknown/forgotten nowadays. These methods belong to boundary-type solution procedures. Below I shall quote a Professor\footnote{His/Her identity will be hidden in order to avoid any kind of diplomatic incidents.} from Laboratoire Jacques-Louis Lions (LJLL) at Paris 6 Pierre and Marie \textsc{Curie} University, who makes interesting comments about the knowledge of \textsc{Trefftz} methods in French Applied Mathematics community (the original language and orthography are conserved):

\begin{quote}

[$\,$\dots] Sinon j'ai moi aussi fait une petite enqu\^ete, \`a la suite de laquelle il appara\^it que personne ne connait \textsc{Trefftz} chez les matheux de Paris 6, hormis \textsc{Nataf} qui en a entendu parler pendant son DEA en m\'eca !!!

\textsc{Trefftz} est un grand oubli\'e car son papier est vraiment extr\^emement int\'eressant, encore  plus si tu te rends compte qu'il produit en fait une estimation \emph{a posteriori} (\c{c}a doit m\^eme \^etre la premi\`ere). En revanche je me suis aussi persuad\'e que les m\'ethodes de \textsc{Trefftz} v\'eg\`etent chez les m\'ecano. [$\,$\dots]

\end{quote}

So, the situation is rather sad. Fortunately this method continued to live in the Mechanical Engineering community who preserved it for future generations.

In the exposition below we shall follow an excellent review paper \cite{Kita1995}. More precisely we describe the so-called \emph{indirect Trefftz method}, which was proposed originally by E.~\textsc{Trefftz} (1926) in \cite{Trefftz1926}. There are also \emph{direct Trefftz methods} proposed some sixty years later in \cite{Cheung1989}. They are much closer to the Boundary Integral Equation Methods (BIEM) and will not be covered here. Interested readers can refer to \cite{Kita1995}.

\begin{figure}
  \centering
  \includegraphics[width=1.0\textwidth]{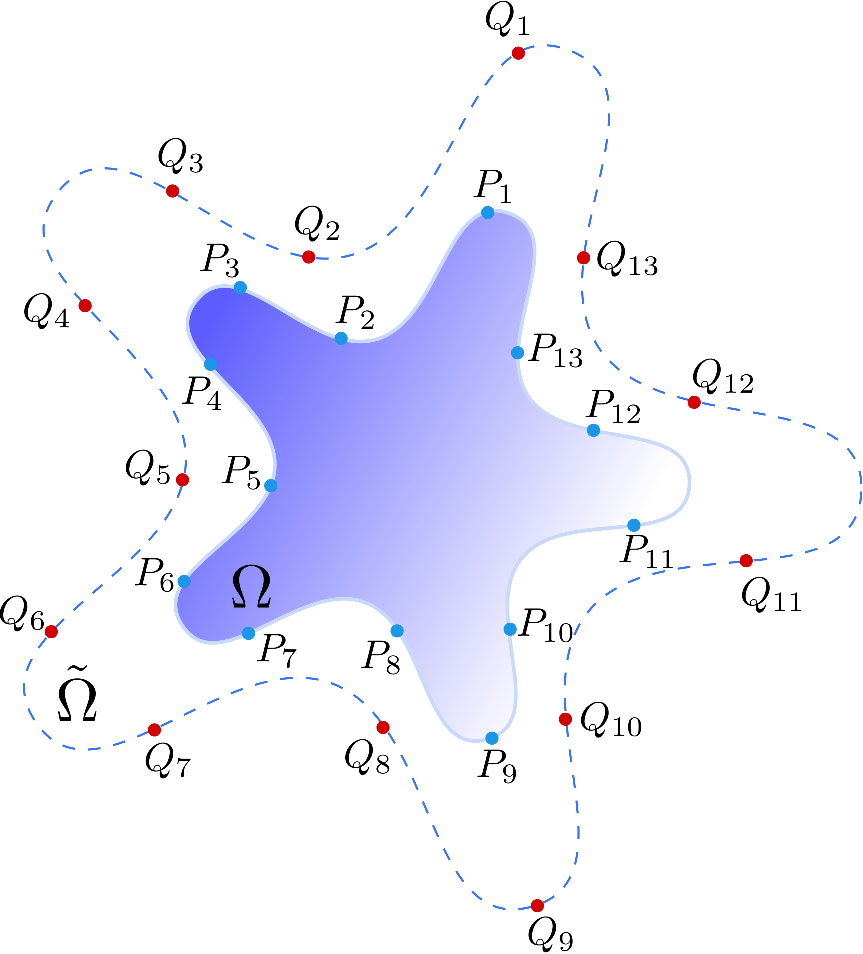}
  \caption{\small\em Collocation Trefftz method using the Green function $\Gr(P,Q)$: $\{P_i\}_{i = 1}^{13}\ \subseteq\ \R^d$ are collocation points and $\{Q_i\}_{i = 1}^{13}\ \subseteq\ \R^d$ are the sources put outside of the domain $\Omega$ in order to avoid singularities. Then, the approximate solution is sought as a linear combination of functions $\bigl\{\Gr(P, Q_i)\bigr\}_{i = 1}^{13}$. The unknown expansion coefficients are found so that to satisfy the boundary conditions exactly in collocation points $\{P_i\}_{i = 1}^{13}$. The outer points $\{Q_i\}_{i = 1}^{13}$ are not needed if we know a $\Tt$-complete set of functions for the operator $\Ll$.}
  \label{fig:imag}
\end{figure}

Consider a (compact) domain $\Omega\ \subseteq\ \R^d$, $d\ \geq\ 2$ and a (scalar) equation on it (see Figure~\ref{fig:imag} for an illustration):
\begin{equation}\label{eq:lapl}
  \Ll u\ =\ 0\,.
\end{equation}
You can think, for example, that $\Ll\ =\ -\grad^2\ =\ -\sum_{i\, =\, 1}^{d}\, \pd{^2}{x_i^2}$ is the \textsc{Laplace}\footnote{Pierre--Simon de \textsc{Laplace} (1749 -- 1827) is a French Mathematician whose works were greatly disregarded during his times. The understanding of their importance came much later (il vaut mieux tard que jamais \smiley).} operator. Equation~\eqref{eq:lapl} has to be completed by appropriate boundary conditions:
\begin{equation}\label{eq:bound}
  \B\, u\ =\ u^\circ\,, \qquad \x\ \in\ \partial\,\Omega\,,
\end{equation}
where $u^\circ$ is the boundary data (solution value or its flux) and $\B$ is an operator depending on the type of boundary conditions in use. For example, in the case of \textsc{Laplace} equation the following choices are popular:
\begin{description}
  \item[Dirichlet] $\B\ \eqdef\ \Id$ (Identity operator, \ie~$\Id u\ \equiv\ u$); in this case the solution value is prescribed on the boundary.
  \item[Neumann] $\B\ \eqdef\ \pd{}{n}\ \equiv\ \n\scal\grad\ =\ \sum_{i\, =\, 1}^{d}\, n_i\cdot\pd{}{x_i}$ (Normal derivative, $\n$ being the exterior normal to $\partial\,\Omega$); physically it corresponds to the prescribed flux through the boundary.
  \item[Robin] $\B\ \eqdef\ \alpha\,\Id\ +\ \beta\;\pd{}{n}$, where $\alpha,\, \beta\ \in\ \R$ are some parameters; this case is a mixture of two previous situations. It arises in some problems as well (see \eg \cite{Chhay2015} for the heat conduction problem in thin liquid films).
\end{description}
Notice, that the boundary $\partial\,\Omega$ can be divided in some sub-domains where a different boundary condition is imposed:
\begin{equation*}
  \partial\,\Omega\ =\ \Gamma_1\ \cup\ \Gamma_2\ \cup\ \ldots\ \cup\ \Gamma_n\,
  \qquad \mu_{d}\,\Bigl(\Gamma_i\ \cap\ \Gamma_j\Bigr)\ =\ 0\,, \quad 1\ \leq\ i\ <\ j\ \leq\ n\,,
\end{equation*}
where $\mu_d$ is the \textsc{Lebesgue} measure in $\R^d$.

In the indirect \textsc{Trefftz} method the numerical solution is sought as a linear combination of $\Tt$-complete functions $\bigl\{\phi_k(\x)\bigr\}_{k\,=\,1}^{n}$ \cite{Herrera1984}, which satisfy exactly the governing equation \eqref{eq:lapl}:
\begin{equation}\label{eq:approx}
  u_n(\x)\ =\ \sum_{k\, =\, 1}^{n}\, v_k\, \phi_k(\x)\,.
\end{equation}
For example, for the \textsc{Laplace} equation in $\R^2$ the $\Tt$-complete set of functions is \cite{Herrera1982} $\bigl\{r^{\,k}\,\ue^{\,\ui k \theta}\bigr\}_{k = 0}^{\infty}$, where $(r,\,\theta)$ are polar coordinates\footnote{$r\ =\ \sqrt{x^2\ +\ y^2}$, $\quad\theta\ =\ \arctan\,\dfrac{y}{x}\,,\quad$ with $\arctan\,(\pm\infty)\ =\ \pm\dfrac{\pi}{2}\,$.} on plane. The coefficients $\{v_k\}_{k = 1}^{n}$ are to be determined. To achieve this goal the approximate solution \eqref{eq:approx} is substituted into the boundary conditions \eqref{eq:bound} to form the residual:
\begin{equation*}
  \Rr[u_n]\ =\ \B\,u_n\ -\ u^\circ\,.
\end{equation*}
If the residual $\Rr[u_n]$ is equal identically to zero on the boundary $\partial\,\Omega$, then we found the exact solution (we are lucky!). Otherwise (we are unlucky and in Numerical Analysis it happens more often), the residual has to be minimized. Very often the boundary operator $\B$ is linear and we can write:
\begin{equation*}
  \Rr[u_n]\ =\ \sum_{k\, =\, 1}^{n}\, v_k\, \B\,\phi_k(\x)\ -\ u^\circ\,.
\end{equation*}
So, we can apply now the collocation, \textsc{Galerkin} or least square methods to determine the coefficients $\{v_k\}_{k\, =\, 1}^{n}$. These procedures were explained above.

\begin{remark}
Suppose that a $\Tt$-complete set of functions for equation \eqref{eq:lapl} is unknown. However, there is a way to overcome this difficulty if we know the \textsc{Green}\footnote{George \textsc{Green} (1793 -- 1841), a British Mathematician who made great contributions to the Mathematical Physics and PDEs.} function\footnote{The \textsc{Green} function is named after George~\textsc{Green} (see the footnote above), but the notion of this function can be already found first in the works of \textsc{Laplace}, then in the works of \textsc{Poisson}. As I said, the works of \textsc{Laplace} were essentially disregarded by his ``colleagues''. \textsc{Poisson} was luckier in this respect. He is also known for (mis)using his administrative resource in order to delay (just for a couple of decades, nothing serious \smiley) the publication of his competitors, \eg the young (at that time) \textsc{Cauchy}.} $\Gr(\x;\,Q)$, $\x \in\ \Omega$ for our equation \eqref{eq:lapl}, \ie
\begin{equation*}
  \Ll\,\Gr\ =\ \delta\bigl(\x\ -\ Q\bigr)\,,\qquad \x\ \in\ \Omega\,, \qquad Q\ \in\ \R^d\,.
\end{equation*}
where $\delta(\x)$ is the singular \textsc{Dirac} measure. Notice that point $Q$ can be inside or outside of domain $\Omega$, where the problem is defined. Then, we can seek for an approximate solution $u_n(\x)$ as a linear combination of \textsc{Green} functions:
\begin{equation*}
  u_n(\x)\ =\ \sum_{k\, =\, 1}^{n}\, v_k\, \Gr_k(\x;\, Q_k)\,,
\end{equation*}
where the points $\bigl\{Q_k\bigr\}_{k = 1}^{n}$ are distributed outside of the computational domain $\Omega$ in order to avoid the singularities. This method is schematically illustrated in Figure~\ref{fig:imag}. The accuracy of this approach depends naturally on the distribution of points $\bigl\{Q_k\bigr\}_{k = 1}^{n}$. To Author's knowledge there exist no theoretical indications for the optimal distribution (as it is the case of \textsc{Tchebyshev} nodes on a segment). So, it remains mainly an experimental area of the research.
\end{remark}

\bigskip
\paragraph*{Conclusions.}

In the indirect \textsc{Trefftz} methods (historically, the original one) the problem solution is sought as a linear combination of the functions, which satisfy the governing equation identically. Then, the unknown coefficients are chosen so that the approximate solution satisfies the boundary condition(s) as well. It can be done by means of collocation, \textsc{Galerkin} or least square procedures. The main advantage of \textsc{Trefftz} methods is that it allows to reduce problem's dimension by one (\ie~3D $\rightsquigarrow$ 2D and 2D $\rightsquigarrow$ 1D), since the residual $\Rr[u_n]$ is minimized on the domain boundary $\partial\,\Omega\ \subseteq\ \R^{d-1}$. Readers who are interested in the application of \textsc{Trefftz} methods to their problems should refer to \cite{Kita1995} and references therein.


\section{A brief history of diffusion in Physics}
\label{app:hist}

Since the main focus of the PhD school is set on the diffusion processes (molecular diffusion, heat and moisture conduction through the walls, \etc), it is desirable to explain how this research started and why the diffusion is generally modeled by \emph{parabolic PDEs} \cite{Evans2010}. The historic part of this Appendix is partially based on \cite{Philibert2005}.

Firs of all, let us make a general, perhaps surprising, remark, which is based on the analysis of historical investigations: understanding the microscopic world is not compulsory to propose a reliable macroscopic law. Indeed, \textsc{Fourier} did not know anything about the nature of heat, \textsc{Ohm}\footnote{Georg Simon~\textsc{Ohm} (1789 -- 1854) was a German Physicist.} about the nature of electricity, \textsc{Fick}\footnote{Adolf Eugen~\textsc{Fick} (1829 -- 1901) was a German physician and physiologist.} about salt solutions and \textsc{Darcy}\footnote{Henry~\textsc{Darcy} (1803 - 1858) was a French engineer in Hydraulics.} about the structure of porosity and water therein. In cases of the diffusion, the bridge between microscopic and macroscopic worlds was built by A.~\textsc{Einstein}\footnote{Albert~\textsc{Einstein} (1879 -- 1955) was a German theoretical Physicist. This personality does not need to be introduced.} in 1905. Namely, he expressed a macroscopic quantity --- the diffusion coefficient --- in terms of microscopic data (this result will be given below). In general, 1905 was \emph{Annus Mirabilis}\footnote{\emph{Annus Mirabilis} comes from Latin and stands for ``extraordinary year'' in English or ``Wunderjahr'' in German.} for A.~\textsc{Einstein}. During this year he published four papers in \emph{Annalen der Physik} whose aftermath was prodigious.

The first scientific study of diffusion was performed by a Scottish chemist Thomas~\textsc{Graham} (1805 -- 1869). His research on diffusion was conducted between 1828 and 1833. Here we quote his first paper:

\begin{quote}
[$\,$\dots] the experimental information we possess on the subject amounts to little more than the well established fact, that gases of different nature, when brought into contact, do not arrange themselves according to their density, the heaviest undermost, and the lighter uppermost, but they spontaneously diffuse, mutually and equally, through each other, and so remain in the intimate state of mixture for any length of time.
\end{quote}

In 1867 James~\textsc{Maxwell}\footnote{James Clerk~\textsc{Maxwell} (1831 -- 1879) was a Scottish Physicist, famous for \textsc{Maxwell} equations.}  estimated the diffusion coefficient of $\mathrm{CO}_2$ in the air using the results (measurements) of \textsc{Graham}. The resulting number was obtained within 5\% accuracy. It is pretty impressive.

A.~\textsc{Fick}\footnote{Adolf~\textsc{Fick}  was the author of the first treatise of \emph{Die medizinische Physik} (Medical Physics) (1856) where he discussed biophysical problems, such as the mixing of air in the lungs, the work of the heart, the heat economy of the human body, the mechanics of muscular contraction, the hydrodynamics of blood circulation, \etc} hold a chair of physiology in W\"urzburg for 31 years. His main contributions to Physics were made during a few years around 1855. When he was 26 years old he published a paper on diffusion establishing the now classical \textsc{Fick}'s diffusion law. \textsc{Fick} did not realize that dissolution and diffusion processes result from the movement of separate entities of salt and water. He deduced the quantitative law proceeding in analogy with the work of J.~\textsc{Fourier} \cite{Fourier1822} who modeled the heat conduction:

\begin{quote}
[\,\dots] It was quite natural to suppose that this law for diffusion of a salt in its solvent must be identical with that according to which the diffusion of heat in a conducting body takes place; upon this law \textsc{Fourier} founded his celebrated theory of heat, and it is the same that \textsc{Ohm} applied [\,\dots] to the conduction of electricity [\,\dots] according to this law, the transfer of salt and water occurring in a unit of time between two elements of space filled with two different solutions of the same salt, must be, \emph{ceteris partibus}, directly proportional to the difference of concentrations, and inversely proportional to the distance of the elements from one another.
\end{quote}

Going along this analogy, \textsc{Fick} assumed that the flux of matter is proportional to its concentration gradient with a proportionality factor $\kappa$, which he called ``\emph{a constant dependent upon the nature of the substances}''. Actually, \textsc{Fick} made an error in (minus) sign which introduces the anti-diffusion and leads to an ill-posed problem. This error was not a source of difficulty for \textsc{Fick} since he analyzed only steady states in accordance with available experimental conditions at that time.

The fundamental article \cite{Einstein1905} devoted to the study of Brownian\footnote{Robert~\textsc{Brown} (1773 -- 1858) was a Scottish Botanist.} motion and entitled ``On the Motion of Small Particles Suspended in a Stationary Liquid, as Required by the Molecular Kinetic Theory of Heat'' by \textsc{Einstein} was published on 18\up{th} July 1905 in \emph{Annalen der Physik}. By the way, it is the most cited \textsc{Einstein}'s paper among four works published during his \emph{Annus Mirabilis}. This manuscript is of capital interest to us as well. Let us quote a paragraph from \cite{Einstein1905}:

\begin{quote}
In this paper it will be shown that, according to the molecular kinetic theory of heat, bodies of a microscopically visible size suspended in liquids must, as a result of thermal molecular motions, perform motions of such magnitudes that they can be easily observed with a microscope. It is possible that the motions to be discussed here are identical with so-called Brownian molecular motion; however, the data available to me on the latter are so imprecise that I could not form a judgment on the question [$\,$\dots]
\end{quote}

\textsc{Einstein} was the first to understand that the main quantity of interest is the mean square displacement $\langle X^2(t)\rangle$ and not the average velocity of particles. Taking into account the discontinuous nature of particle trajectories, the velocity is meaningless.

Before publication of \cite{Einstein1905} atoms in Physics were considered as a useful, but purely theoretical concept. Their reality was seriously debated. W.~\textsc{Ostwald}\footnote{Friedrich Wilhelm~\textsc{Ostwald} (1853 -- 1932) was a Latvian Chemist. He received a Nobel prize in Chemistry in 1909 for his works on catalysis.}, one of the leaders of the anti-atom school, later told A.~\textsc{Sommerfeld}\footnote{Arnold Johannes Wilhelm~\textsc{Sommerfeld} (1868 -- 1951) was a German theoretical Physicist. He served as the PhD advisor for several Nobel prize winners.} that he had been convinced of the existence of atoms by \textsc{Einstein}'s complete explanation of Brownian motion.

Without knowing it, \textsc{Einstein} in  answered a question posed the same year by K.~\textsc{Pearson}\footnote{Karl~\textsc{Pearson} (1857 -- 1936) was an English Statistician.} in a Letter published in \textsc{Nature} \cite{Pearson1905}:

\begin{quote}
Can any of your readers refer me to a work wherein I should find a solution of the following problem, or failing the knowledge of any existing solution provide me with an original one? I should be extremely grateful for aid in the matter.

A man starts from a point $O$ and walks $\ell$ yards in a straight line; he then turns through any angle whatever and walks another $\ell$ yards in a second straight line. He repeats this process $n$ times. I require the probability that after these $n$ stretches he is at a distance between $r$ and $r\ +\ \ud r$ from his starting point, $O$.

The problem is one of considerable interest, but I have only succeeded in obtaining an integrated solution for two stretches. I think, however, that a solution ought to be found, if only in the form of a series in powers of $\dfrac{1}{n}$, when $n$ is large.
\end{quote}
The expression ``\emph{random walk}'' was probably coined in \textsc{Pearson}'s Letter \cite{Pearson1905}.

Let us follow \textsc{Einstein}'s study of the Brownian motion. Consider a long thin tube filled with water (it allows us to consider only one spatial dimension). At the initial time $t\ =\ 0$ we inject a unit amount of ink at $x\ =\ 0$. Let $u(x,\,t)$ denote the density of ink particles in location $x\ \in\ \R$ and at time $t\ \geq\ 0$. So, initially we have
\begin{equation}\label{eq:initcond}
  u(x,\,0)\ =\ \delta(x)\,,
\end{equation}
where $\delta(x)$ is the \textsc{Dirac} distribution centered at zero.

Then, suppose that the probability density of the event that an ink particle moves from $x$ to $x + \Delta x$ in a small time $\Delta t$ is $\rho(\Delta x,\, \Delta t)$. Then, we have
\begin{equation*}
  u\,(x,\, t\ +\ \Delta t)\ =\ \int_{-\infty}^{+\infty}\, u(x - \Delta x,\, t)\, \rho(\Delta x,\, \Delta t)\; \ud\bigl(\Delta x\bigr)\,.
\end{equation*}
Assuming that solution $u(x,\,t)$ is smooth, we can apply the \textsc{Taylor}\footnote{Brook~\textsc{Taylor} (1685 -- 1731) was an English Mathematician.} formula:
\begin{equation}\label{eq:hh}
  u\,(x,\, t\ +\ \Delta t)\ =\ \int_{-\infty}^{+\infty}\,\Bigl(u\ -\ u_x\cdot\Delta x\ +\ \frac{1}{2}\, u_{xx}\cdot \bigl(\Delta x\bigr)^2\ +\ \cdots\Bigr)\, \rho(\Delta x,\, \Delta t)\; \ud\bigl(\Delta x\bigr)\,.
\end{equation}
Now let us recall that $\rho$ is a Probability Density Function (PDF). Thus,
\begin{description}
  \item[Normalisation] $\int_{-\infty}^{+\infty}\, \rho(\Delta x,\, \Delta t)\; \ud\bigl(\Delta x\bigr)\ =\ 1$
  \item[Symmetry] we can assume that the diffusion process is symmetric in space, \ie
  \begin{equation*}
    \rho(\Delta x,\, \Delta t)\ =\ \rho(-\Delta x,\, \Delta t)\,, \qquad \forall \Delta x,\ \Delta t\ \geq\ 0\,.
  \end{equation*}
  Consequently,
  \begin{equation*}
    \int_{-\infty}^{+\infty}\,\Delta x\, \rho(\Delta x,\, \Delta t)\; \ud\bigl(\Delta x\bigr)\ =\ 0\,.
  \end{equation*}
  \item[Variance] We can also assume that the variance if finite and \textsc{Einstein} assumed additionally that the variance is \emph{linear} in $\Delta t$:
  \begin{equation*}
    \int_{-\infty}^{+\infty}\,\bigl(\Delta x\bigr)^2\, \rho(\Delta x,\, \Delta t)\; \ud\bigl(\Delta x\bigr)\ =\ \D \Delta t\,,
  \end{equation*}
  where $\D\ >\ 0$ is the so-called \emph{diffusion constant}.
\end{description}
By incorporating these results into equation \eqref{eq:hh} and rearranging the terms, it becomes
\begin{equation*}
  \frac{u(x,\,t + \Delta t)\ -\ u(x,\,t)}{\Delta t}\ =\ \frac{1}{2}\,\D\, u_{xx}\,(x,\,t)\ +\ \cdots\,.
\end{equation*}
Taking the limit $\Delta t\ \to\ 0$ we obtain straightforwardly
\begin{equation}\label{eq:par}
  u_t\ =\ \frac{1}{2}\,\D\, u_{xx}\,.
\end{equation}
The Initial Value Problem \eqref{eq:initcond} for this linear parabolic equation \eqref{eq:par} can be solved exactly:
\begin{equation*}
  u(x,\,t)\ =\ \frac{1}{\sqrt{2\pi \D t}}\, \ue^{-\frac{x^2}{2\D t}}\,.
\end{equation*}
The last solution is also known as the \textsc{Green} function.

However, the main result of \textsc{Einstein}'s paper \cite{Einstein1905} is the following formula:
\begin{equation}\label{eq:einstein}
  \D\ =\ \frac{RT}{f N_a}\,,
\end{equation}
where
\begin{description}
  \item[$R$] the ideal gas constant, \ie~$R\ \approx\ 8.3144598$ $\dfrac{\mathsf{J}}{\mathsf{K}\cdot \mathsf{mol}}$
  \item[$T$] the absolute temperature
  \item[$f$] the friction coefficient
  \item[$N_a$] \textsc{Avogadro}'s\footnote{Lorenzo Romano Amedeo Carlo \textsc{Avogadro} di Quaregna e di Cerreto (1776 -- 1856) was an Italian scientist.} number, \ie~$N_a\ \approx\ 6.022140857$ $\mathsf{mol}^{-1}$
\end{description}
Formula \eqref{eq:einstein} along with the observation of the Brownian motion enabled J.~\textsc{Perrin}\footnote{Jean Baptiste~\textsc{Perrin} (1870 -- 1942) was a French Physicist who was honoured with the Nobel prize in Physics in 1926 for the confirmation of the atomic nature of matter (through the observation of Brownian motion).} to produce the first historical estimation of \textsc{Avogadro}'s constant $N_a\ \approx\ 6$ $\mathsf{mol}^{-1}$.

\begin{exo}
Read a fascinating paper by A.~\textsc{Einstein} on the Hydrodynamics of tea leaves \cite{Einstein1926}. By the way, the Author attests that Tea is very stimulating for intellectual activities. \smiley
\end{exo}


\section{Monte--Carlo approach to the diffusion simulation}
\label{app:mc}

From considering the history and physical modelling of Brownian motion we naturally come to the numerical simulation of parabolic PDEs using stochastic processes. It falls into the large class of Monte--Carlo and quasi-Monte--Carlo methods \cite{Caflisch1998}. They are based on the so-called \emph{Law of large numbers}, which can be informally stated\footnote{We choose this particular form, since it is suitable for our exposition.} as
\begin{equation*}
  \lim_{N \to +\infty} \frac{1}{N}\sum_{n\, =\, 1}^{N} g(\xi_n)\ =\ \E\,[g(\xi)]\,,
\end{equation*}
where $\xi_n$ are independent, identically distributed random variables and $g(\cdot)$ is a real-valued continuous function. Here $\E\,[\cdot]$ stands for the mathematical expectation and the convergence is understood in the sense \emph{convergence in probability distribution}\footnote{A sequence of random variables $\{\xi_n\}_{n = 1}^{\infty}$ \emph{converges in probability distribution} towards the random variable $\xi$ if for $\forall\, \eps\ >\ 0$
\begin{equation*}
  \lim_{n \to\ \infty} \P\bigl\{\abs{\xi_n\ -\ \xi}\ \geq\ \eps\bigr\}\ =\ 0\,.
\end{equation*}
This fact can be denoted as $\xi_n\ \stackrel{\P}{\rightsquigarrow}\ \xi$.}. The proof is based on the \textsc{Tchebyshev} inequality in Probabilities, but we do not enter into these details. The main advantages of Monte--Carlo methods are
\begin{itemize}
  \item Simplicity of implementation
  \item Independent of the problem dimension (they do not suffer of the curse of dimensionality)
\end{itemize}
On the other hand, the convergence rate is given by the Central Limit Theorem \cite{Jaynes2003} and it is rather slow, \ie~$\O(N^{\,-1/2})$, even if some acceleration is possible thanks to some adaptive procedures \cite{Lapeyre2011} (such as low-discrepancy sequences, variance reduction and multi-level methods). The large deviation theory guarantees that the probability of falling out of a fixed tolerance interval decays exponentially fast.

To the Author knowledge, the method we are going to describe below was first proposed by R.~\textsc{Feynman}\footnote{Richard Phillips~\textsc{Feynman} (1918 -- 1988) was an American theoretical Physicist. Nobel Prize in Physics (1965) and Author's hero. Please, do not hesitate to read any of his books!} in order to solve \emph{numerically} the linear \textsc{Schr\"odinger}\footnote{Erwin Rudolf Josef Alexander~\textsc{Schr\"odinger} (1887 -- 1961) was an Austrian theoretical Physicist. Nobel prize in Physics (1933) for the formulation of what is known now as the \textsc{Schr\"odinger} equation.} equation in 1940s. In fact, he noticed that the \textsc{Schr\"odinger} equation can be solved by a kind of average over trajectories. This observation led him to a far-reaching reformulation of the quantum theory in terms of \emph{path integrals}. Upon learning \textsc{Feynman}'s ideas, M.~\textsc{Kac}\footnote{Mark~\textsc{Kac} (1914 -- 1984) was a Polish mathematician who studied in Lviv University, Ukraine and immigrated later to USA.} (\textsc{Feynman}'s colleague at Cornell University) understood that a similar method can work for the heat equation as well. Later it became \textsc{Feynman}--\textsc{Kac} method. Unfortunately, this method is not implemented in any commercial software (for PDEs, the financial industry is using these methods since many decades) and it is not described in classical textbooks on PDEs and in the books on numerical methods for PDEs.

Consider first the following simple heat equation:
\begin{equation}\label{eq:heat1}
  u_t\ =\ \frac{1}{2}\; u_{xx}\ -\ V(x)\cdot u\,,
\end{equation}
where $V(x)$ is a function representing the amount of external cooling (if $V(x) \geq 0$, and external heating if $V(x) < 0$, but we do not consider this case) at point $x$. Then, we have the following
\begin{theorem}[\textsc{Feynman}--\textsc{Kac} formula]\label{thm:fk}
Let $V(x)$ be a non-negative continuous function and let $u_0(x)$ be bounded and continuous. Suppose that $u(x,\,t)$ is a bounded function that satisfies equation \eqref{eq:heat1} along with the initial condition
\begin{equation}\label{eq:init1}
  u(x,\,0)\ =\ u_0(x)\,,
\end{equation}
then,
\begin{equation}\label{eq:fk}
  u(x,\,t)\ =\ \E\,\biggl[\,\exp\Bigl\{-\int_0^t V\bigl(\W_s\bigr)\;\ud s\Bigr\}\; u_0\bigl(\W_t\bigr)\,\biggr]\,,
\end{equation}
where $\{\W_t\}_{t\,\geq\,0}$ is a Brownian motion starting at $x$.
\end{theorem}
We have to define rigorously also what the Brownian motion is:
\begin{deff}[Brownian motion]
A real-valued stochastic process $t\ \mapsto\ \W(t)$ (or a random curve) is called a \emph{Brownian motion} (or \emph{Wiener process}) if
\begin{enumerate}
  \item Almost surely $\W(0)\ =\ 0$
  \item $\W(t)\ -\ \W(s)\ \sim\ \Nn(0,\, t\ -\ s)\ \equiv\ \sqrt{t\ -\ s}\,\Nn(0,\, 1)$, $\quad\forall t\ \geq\ s\ \geq\ 0$
  \item For any instances of time $v\ >\ u\ >\ t\ >\ s\ \geq\ 0$ the increments $\W(v) - \W(u)$ and $\W(t) - \W(s)$ are independent random variables.
\end{enumerate}
\end{deff}
In particular, from point (2) by choosing $s = 0$ it follows directly that
\begin{equation*}
  \E\bigl[\W(t)\bigr]\ =\ 0, \qquad
  \E\bigl[\W^{\,2}(t)\bigr]\ =\ t\,.
\end{equation*}
More informally we can say that the Brownian motion is a continuous random curve with the \emph{largest} possible amount of randomness.

The Theorem above can be proved even under more general assumptions, but the formulation given above suffices for most practical situations. For instance, the functions $V(x)$ and $u_0(x)$ may have isolated discontinuities and the \textsc{Feynman}--\textsc{Kac} formula will be still valid. Another interesting corollary of \textsc{Feynman}--\textsc{Kac}'s formula is the uniqueness of solutions to the Initial Value Problem (IVP) \eqref{eq:heat1}:
\begin{corollary}
Under the assumptions of Theorem~\ref{thm:fk}, there is at most one solution to the heat equation \eqref{eq:heat1}, which satisfies the initial condition \eqref{eq:init1}. Namely, it is given by \textsc{Feynman}--\textsc{Kac} formula \eqref{eq:fk}.
\end{corollary}

The main (computational) advantage of \textsc{Feynman}--\textsc{Kac} formula is that it can be straightforwardly generalized to the arbitrary dimension:
\begin{theorem}[\textsc{Feynman}--\textsc{Kac} formula in $d$ dimensions]
Let $V:\ \R^d\ \mapsto\ [0, +\infty)$ and $u_0:\ \R^d\ \mapsto\ \R$ be continuous functions with bounded initial condition $u_0(\x)$. Suppose that $u(\x,\,t)$ is also a bounded function, which satisfies the following Partial Differential Equation (PDE):
\begin{equation*}
  u_t\ =\ \frac{1}{2}\;\sum_{i\, =\, 1}^{d}\,\pd{^2 u}{x_i^{\,2}}\ -\ V(\x)\cdot u\,,
\end{equation*}
and the initial condition
\begin{equation*}
  u(\x,\,0)\ =\ u_0(\x)\,.
\end{equation*}
Then,
\begin{equation*}
  u(\x,\,t)\ =\ \E\,\biggl[\,\exp\Bigl\{-\int_0^t V\bigl(\Wb_s\bigr)\;\ud s\Bigr\}\; u_0\bigl(\Wb_t\bigr)\,\biggr]\,,
\end{equation*}
where $\{\Wb_t\}_{t\,\geq\,0}$ is a $d$-dimensional Brownian motion starting at $\x$.
\end{theorem}

We give here another useful generalization of the \textsc{Feynman}--\textsc{Kac} method. For the sake of simplicity we return to the one-dimensional case. Consider a stochastic process $X_t$, which satisfies the following \emph{Stochastic Differential Equation} (SDE):
\begin{equation}\label{eq:sde}
  \ud X_t\ =\ \alpha(X_t)\,\ud t\ +\ \sigma(X_t)\;\ud\W_t\,, \qquad X_0\ =\ x\,,
\end{equation}
where $\alpha(x)$ is the local drift and $\sigma(x)$ is called the \emph{local volatility} in financial applications. We assume here that $\alpha(x)$ and $\sigma(x)$ are globally Lipschitz continuous and grow linearly in space at most. Then, the function given by formula \eqref{eq:fk} satisfies the generalized diffusion equation
\begin{equation*}
  u_t\ =\ \alpha(x)\,u_x\ +\ \frac{1}{2}\;\sigma^2(x)\,u_{xx}\ -\ V(x)\, u(x)\,,
\end{equation*}
together with the initial condition \eqref{eq:init1}.

So, the resulting numerical algorithm is very simple:
\begin{enumerate}
  \item Generate $N$ trajectories of the Brownian motion $\{\W_s^k\}_{s\, \in\, [0,\,t]}^{1\, \leq\, k\, \leq\, N}$
  \item Compute $N$ solutions to the SDE \eqref{eq:sde} using the \textsc{Euler}--\textsc{Maruyama}\footnote{Gisiro~\textsc{Maruyama} (1916 -- 1986) was a Japanese Mathematician with notable contributions to the theory of stochastic processes.} method \cite{Higham2001}, for example
  \item Compute the solution value using representation \eqref{eq:fk}. The mathematical expectation $\E\,[\cdot]$ is replaced by the simple arithmetic average according to the Monte--Carlo approach.
\end{enumerate}

\begin{remark}
Notice that the step (2) can be omitted if we solve the linear parabolic equation without the drift term and with constant diffusion coefficient. The simplest version of the \textsc{Feynman}--\textsc{Kac} method for the simplest initial value problem
\begin{equation*}
  u_t\ =\ \frac{1}{2}\;u_{xx}\,, \qquad u(x,\,0)\ =\ u_0(x)\,,
\end{equation*}
looks like
\begin{equation*}
  u(x,\,t)\ \approx\ \frac{1}{N}\;\sum_{n\, =\, 1}^{N} u_0(\xi_n)\,, \qquad \forall n:\ \xi_n\ \sim\ \Nn(x,\, \sqrt{t})\,.
\end{equation*}
It is probably the simplest way to estimate solution value in a given point $(x,\,t)$.
\end{remark}

The \textsc{Feynman}--\textsc{Kac} method inherits all advantages (and disadvantages) of Monte--Carlo methods. Comparing to grid-based methods it enjoys also the locality property. Namely, if you want to compute your solution in a given point $(x,\,t)$, you do not need to know the solution in neighbouring sites (locations). So, if you want to compute numerically the solution only in a few specific points, please, consider this method. It can be competitive with more conventional approaches. Moreover, the financial industry already appreciated the power of these methods.

\subsection{Brownian motion generation}
\label{app:brown}

\begin{figure}
  \centering
  \includegraphics[width=0.79\textwidth]{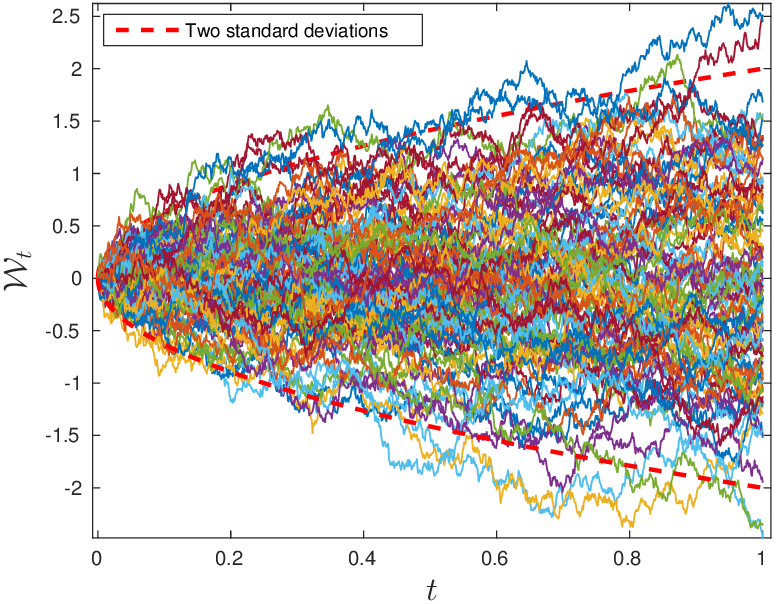}
  \caption{\small\em A collection of sample Brownian paths generated by the code given in Appendix~\ref{app:brown}. Two standard deviations should contain about 99\% of trajectories. On this picture it looks like it is the case.}
  \label{fig:brown1}
\end{figure}

In this Section we provide a simple \textsc{Matlab} code to generate $M$ sample Brownian motion realizations. It can be used as a building block for SDE solvers and practical implementations of the \textsc{Feynman}--\textsc{Kac} method described above:

\begin{lstlisting}
M  = 100;  % number of paths
N  = 1000; % number of steps
T  = 1;    % final simulation time
dt = T/N;  % time step
dW = sqrt(dt)*randn(M, N);
W  = cumsum(dW, 2);
\end{lstlisting}
A sample output result (with precisely the same parameters) of this script is depicted in Figure~\ref{fig:brown1}.

\begin{figure}
  \centering
  \includegraphics[width=0.99\textwidth]{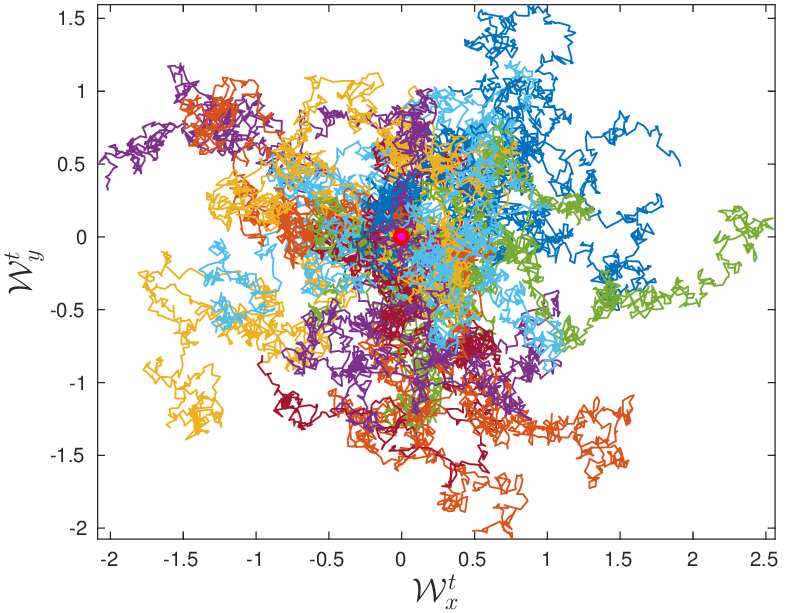}
  \caption{\small\em A few Brownian paths in two spatial dimensions. The red circle depicts the starting point $(0,\,0)$.}
  \label{fig:brown2}
\end{figure}

It is quite straightforward to employ the same \textsc{Matlab} script to generate Brownian paths in $d$ dimensions (just take $M = m\times d$, where $m \in \N$ and then regroup matrix $W$ in $m$ $d$-dimensional trajectories). A few realizations of the Brownian motion in two spatial dimensions are depicted in Figure~\ref{fig:brown2}.

\section{An exact non-periodic solution to the 1D heat equation}

In this Appendix we provide an exact solution (see \cite[Section \textsection 7.2]{Fornberg1996}) to the (dimensionless) heat equation
\begin{equation}\label{eq:heat}
  u_t\ -\ \frac{1}{9\pi^2}\; u_{xx}\ =\ 0\,, \qquad x\ \in\ [0,\, 1]\,, \qquad t\ \in\ \R^+\,,
\end{equation}
subject to the initial condition
\begin{equation*}
  u(x,\,0)\ =\ 0\,, \qquad x\ \in\ [0,\, 1]\,,
\end{equation*}
and the following non-periodic and non-homogeneous boundary conditions:
\begin{equation*}
  u(0,\, t)\ =\ \sin t\,, \qquad u_{\,x}(1,\,t)\ =\ 0\,, \qquad t\ \in\ \R^+\,.
\end{equation*}
The first boundary condition says that the temperature is prescribed at the left boundary and zero heat flux is imposed on the right boundary. So, the unique solution to the problem \eqref{eq:heat} described above for $\forall t\ >\ 0$ is given by\footnote{This solution can be derived using the classical method of separation of variables \cite{Griffiths2015} [or see an excellent publication of Professor (and my favourite colleague and friend) Marguerite~\textsc{Gisclon} \cite{Gisclon1998}, if you are not scared by French].}
\begin{multline}\label{eq:sol}
  u(x,\,t)\ =\ \underbrace{\sinh^{-1}\frac{3\pi}{2}\;\Bigl[\cos\frac{3\pi x}{2}\,\sinh\frac{3 \pi (1 - x)}{2}\, \sin t\ -\ \sin\frac{3 \pi x}{2}\,\cosh\frac{3\pi(1 - x)}{2}\,\cos t\Bigr]}_{u_{\circlearrowleft}(x,\,t)}\\ +\ \underbrace{\frac{72}{\pi}\;\sum_{n\, =\, 1}^{\infty}\,\frac{(2n - 1)\,\ue^{-\frac{(2n - 1)^2}{18}\, t}}{\bigl[9\ +\ 4(n - 2)^2\bigr]\bigl[9\ +\ 4(n + 1)^2\bigr]}\; \sin(n - \half)\pi x}_{u_\Sigma(x,\,t)}\,.
\end{multline}
This analytical solution can be used as the \emph{reference solution} in order to validate your non-periodic numerical codes. Notice that the second term $u_\Sigma(x,\,t)$ (\ie~the infinite series) vanishes uniformly in space, \ie
\begin{equation*}
  u_\Sigma(x,\,t)\ \rightrightarrows\ 0 \quad \mbox{ as }\quad t\ \to\ +\infty\,.
\end{equation*}
This term $u_\Sigma$ is needed to enforce the initial condition. Consequently, the long time behaviour of the solution \eqref{eq:sol} is given by $u_\circlearrowleft(x,\,t)$.

\begin{figure}
  \centering
  \includegraphics[width=0.99\textwidth]{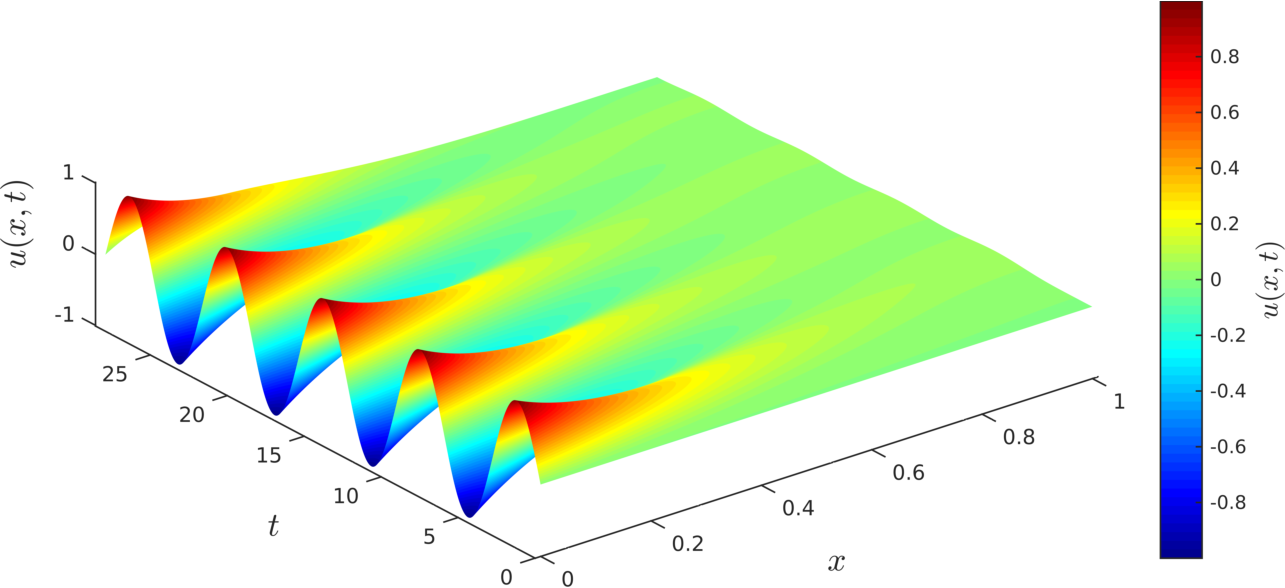}
  \caption{\small\em Space-time plot of the solution \eqref{eq:sol} to the linear heat equation \eqref{eq:heat}. The time window is $t \in\ [0,\, 9\pi]$.}
  \label{fig:heat}
\end{figure}

\begin{figure}
  \centering
  \includegraphics[width=0.99\textwidth]{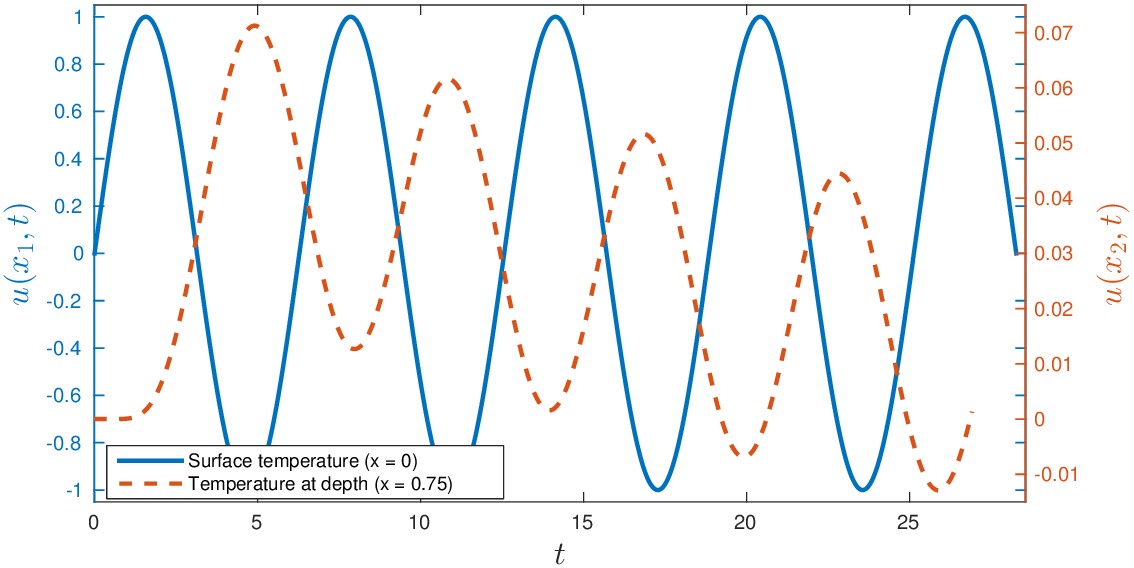}
  \caption{\small\em Solution \eqref{eq:sol} shown at the surface $x = 0$ and at depth $x\ =\ \frac{3}{4}$. The goal is to illustrate the phase shift between two curves.}
  \label{fig:depth}
\end{figure}

\begin{remark}
Equation \eqref{eq:heat} (along with boundary conditions) models the temperature variation in soil under periodic boundary forcing (modeling seasonal temperature variations). Solution \eqref{eq:sol} explains also that for many types of soils, there is a phase shift of seasonal temperature at certain depths (warmest in winter and coolest in summer). See Figure~\ref{fig:depth} for an illustration. It provides a theoretical explanation why a cave works in practice.
\end{remark}

\section{Some popular numerical schemes for ODEs}
\label{app:odes}

The discussion of numerical methods for PDEs cannot be complete if we do not provide some basic techniques for the time marching. Basically, if we follow the Method Of Lines (MOL) \cite{Schiesser1994, Kreiss1992, Reddy1992, Shampine1994}, the PDE is discretized first in space, then the resulting system of coupled ODEs has to be solved numerically. This document would not be complete if we did not provide any indications on how to do it. Consider for simplicity an Initial Value Problem (IVP), which is sometimes also referred to as the \textsc{Cauchy} problem:
\begin{equation*}
  \dot{u}\ =\ f(u)\,, \qquad u(0)\ =\ u_0\,.
\end{equation*}
where the dot over a function denotes the derivative with respect to time, \ie\ $\dot{u}\ \equiv\ \od{u}{t}\,$. There is a number of time marching schemes proposed in the literature. We refer to \cite{Hairer2009, Hairer1996, Hairer2002} as exhaustive references on this topic. Below we provide some most popular (subjectively) schemes.
\begin{description}
  \item[Forward Euler] (explicit, first order accurate)
  \begin{equation*}
    u_{n+1}\ =\ u_n\ +\ \Delta t\,f(u_n)\,.
  \end{equation*}
\end{description}
\begin{description}
  \item[Backward Euler] (implicit, first order accurate)
  \begin{equation*}
    u_{n+1}\ =\ u_n\ +\ \Delta t\,f(u_{n+1})\,.
  \end{equation*}
\end{description}
\begin{description}
  \item[Adams\footnotemark--Bashforth-2]\footnotetext{John Couch \textsc{Adams} (1819 -- 1892), a British Mathematician and Astronomer. He predicted the existence and position of the planet Neptune.} (explicit, second order accurate)
  \begin{equation*}
    u_{n+1}\ =\ u_n\ +\ \Delta t\,\Bigl[\,\frac{3}{2}\;f(u_{n})\ -\ \frac{1}{2}\;f(u_{n-1})\,\Bigr]\,.
  \end{equation*}
\end{description}
\begin{description}
  \item[Adams--Bashforth\footnotemark-3]\footnotetext{Francis \textsc{Bashforth} (1819 -- 1912) is a British Applied Mathematician working in the field of Ballistics. However, his famous numerical scheme was proposed in collaboration with J.~C.~\textsc{Adams} to study the drop formation.} (explicit, third order accurate)
  \begin{equation*}
    u_{n+1}\ =\ u_n\ +\ \Delta t\,\Bigl[\,\frac{23}{12}\;f(u_{n})\ -\ \frac{4}{3}\;f(u_{n-1})\ +\ \frac{5}{12}\;f(u_{n-2})\,\Bigr]\,.
  \end{equation*}
\end{description}
\begin{description}
  \item[Adams--Moulton\footnotemark-1 or the trapezoidal rule]\footnotetext{Forest Ray \textsc{Moulton} (1872 -- 1952) is an American Astronomer. There is a crater on the Moon named after him.} (implicit, second order accurate)
  \begin{equation*}
    u_{n+1}\ =\ u_n\ +\ \frac{1}{2}\;\Delta t\,\Bigl[\,f(u_{n+1})\ +\ f(u_{n})\,\Bigr]\,.
  \end{equation*}
\end{description}
\begin{description}
  \item[Adams--Moulton-2] (implicit, third order accurate)
  \begin{equation*}
    u_{n+1}\ =\ u_n\ +\ \Delta t\,\Bigl[\,\frac{5}{12}\;f(u_{n+1})\ +\ \frac{2}{3}\;f(u_{n})\ -\ \frac{1}{12}\;f(u_{n-1})\,\Bigr]\,.
  \end{equation*}
\end{description}
\begin{description}
  \item[Two-stage Runge--Kutta\footnotemark]\footnotetext{Martin Wilhelm \textsc{Kutta} (1867 -- 1944) is a German Mathematician who worked in the field of Fluid Mechanics and Numerical Analysis.} (explicit, second order accurate)
  \begin{align*}
    k_1\ &=\ \Delta t\,f(u_n)\,, \\
    k_2\ &=\ \Delta t\,f\bigl(u_n\ +\ \alpha\;k_1\bigr)\,, \\
    u_{n+1}\ &=\ u_n\ +\ \bigl(1\ -\ \frac{1}{2\,\alpha}\bigr)\,k_1\ +\ \frac{1}{2\,\alpha}\;k_2\,.
  \end{align*}
  \begin{itemize}
    \item $\alpha\ =\ \frac{1}{2}$: mid-point scheme
    \item $\alpha\ =\ 1$: \textsc{Heun}'s\footnote{Karl~\textsc{Heun} (1859 -- 1929) is a German Mathematician.} method
    \item $\alpha\ =\ \frac{2}{3}$: \textsc{Ralston}'s method
  \end{itemize}
\end{description}
\begin{description}
  \item[Runge--Kutta-4 (RK4)] (explicit, fourth order accurate)
  \begin{align*}
    k_1\ &=\ \Delta t\,f(u_n)\,, \\
    k_2\ &=\ \Delta t\,f\bigl(u_n\ +\ \frac{1}{2}\;k_1\bigr)\,, \\
    k_3\ &=\ \Delta t\,f\bigl(u_n\ +\ \frac{1}{2}\;k_2\bigr)\,, \\
    k_4\ &=\ \Delta t\,f(u_n\ +\ k_3)\,, \\
    u_{n+1}\ &=\ u_n\ +\ \frac{1}{6}\;\bigl[\,k_1\ +\ 2\,k_2\ +\ 2\,k_3\ +\ k_4\,\bigr]\,.
  \end{align*}
\end{description}
We do not discuss here the questions related to the stability of these schemes. This topic is of uttermost importance but out of scope of these Lecture notes. Moreover, we skip also the adaptive embedded \textsc{Runge}--\textsc{Kutta} schemes which can choose the time step to meet some prescribed accuracy requirements \cite{Dormand1980}.

\begin{figure}
  \centering
  \includegraphics[width=0.79\textwidth]{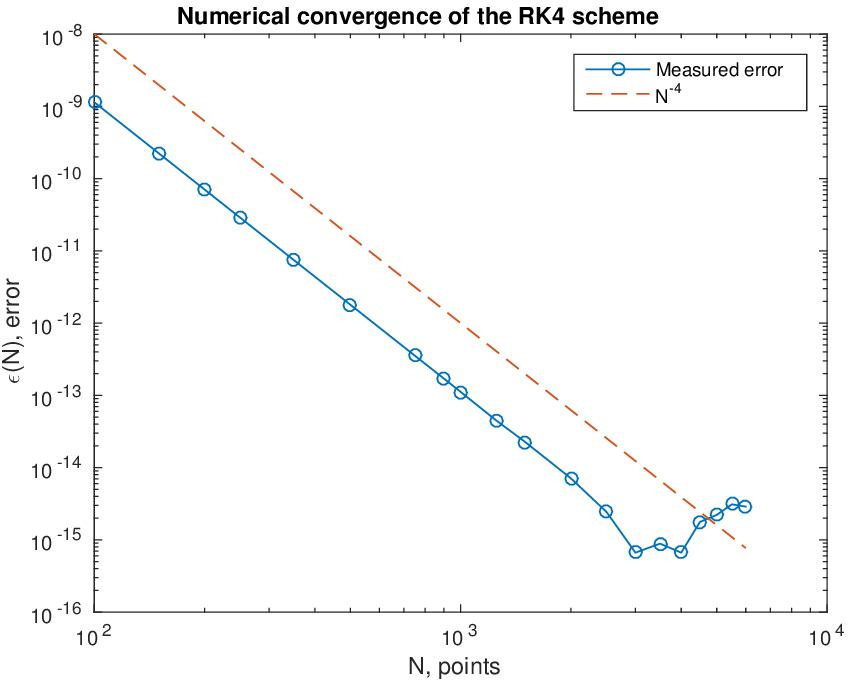}
  \caption{\small\em Numerical convergence of the 4\up{th} order \textsc{Runge}--\textsc{Kutta} scheme when applied to the logistic equation.}
  \label{fig:rk4}
\end{figure}

In order to illustrate the usage of RK4 scheme (sometimes called \emph{the} \textsc{Runge}--\textsc{Kutta} scheme) we take a simple nonlinear logistic equation: 
\begin{equation}\label{eq:logistic}
  \dot{u}\ =\ u\cdot(1\ -\ u)\,, \qquad u(0)\ =\ 2\,.
\end{equation}
It is not difficult to check that the exact solution to this IVP is
\begin{equation*}
  u(t)\ =\ \frac{2}{2\ -\ \ue^{-t}}\,.
\end{equation*}
The \textsc{Matlab} code used to study numerically the convergence of the RK4 scheme on the nonlinear equation \eqref{eq:logistic} is given below:
\newpage
\begin{lstlisting}
T   = 2.0; % final simulation time
uex = 2/(2 - exp(-T));
rhs = @(u) u*(1 - u); % RHS of the logistic equation

% list of successfully refined grids:
NN  = [100; 150; 200; 250; 350; 500; 750; 900; 1000; 1250; 1500; 
       2000; 2500; 3000; 3500; 4000; 4500; 5000; 5500; 6000];

Err = zeros(size(NN));
for n = 1:length(NN)
    N  = NN(n); % number of time steps
    dt = T/N;   % one time step
    u  = 2.0;   % initial condition on the solution
    for j = 1:N
        k1 = dt*rhs(u);
        k2 = dt*rhs(u + 0.5*k1);
        k3 = dt*rhs(u + 0.5*k2);
        k4 = dt*rhs(u + k3);
        u  = u + (k1 + 2*k2 + 2*k3 + k4)/6;
    end % for j
    Err(n) = abs(u - uex);
end % for n
\end{lstlisting}
The numerical result is shown in Figure~\ref{fig:rk4}. One can clearly observe the 4\up{th} order convergence of the RK4 scheme applied to a nonlinear example. The convergence is broken only by the rounding effects inherent to the floating point arithmetics.

\subsection{Existence and unicity of solutions}

Normally, before attempting to solve a differential equation, one has to be sure that the mathematical problem is well-posed. The mathematical notion of the well-posedness was proposed by J.~\textsc{Hadamard}\footnote{Jacques~\textsc{Hadamard} (1865 -- 1963), a French Mathematician who made seminal contributions to several fields of Mathematics including the Number Theory, Complex Analysis and the theory of PDEs.} \cite{Hadamard1902} and it includes three points to be checked:
\begin{enumerate}
  \item Existence
  \item Uniqueness
  \item Continuous dependence on the initial condition and other problem parameters
\end{enumerate}
The last point is usually more difficult to be proven theoretically. However, for 3D \textsc{Navier}\footnote{Claude Louis Marie Henri \textsc{Navier} (1785 -- 1836) is a French Engineer (Corps des Ponts et Chauss\'ees) and Physicist who made contributions to Mechanics.}--\textsc{Stokes}\footnote{George Gabriel \textsc{Stokes} (1819 -- 1903) was an Irish Physicist and Mathematician.} equations even the global existence and smoothness of solutions (\ie the first point) poses already a Millennium problem:
\begin{center}
  \href{http://www.claymath.org/millennium-problems/navier%E2%80%93stokes-equation}{http://www.claymath.org/millennium-problems/navier-stokes-equation}
\end{center}
Below we give some theoretical results which address the well-posedness conditions for an Initial Value Problem (IVP) for a scalar equation:
\begin{equation}\label{eq:ODE}
  \dot{u}\ =\ f(t,\,u)\,, \qquad u(t_0)\ =\ u_0\,.
\end{equation}
The generalization to systems of Ordinary Differential Equations (ODEs) is straightforward.

\begin{theorem}[Existence]
If function $f(t,\,u)$ is continuous and bounded in a domain $(t,\,u)\ \in\ \D\ \subseteq \R_t\times\R_u$, then through every point $(t_0,\,u_0)\ \in\ \D$ passes at least one integral curve of equation \eqref{eq:ODE}.
\end{theorem}

\begin{theorem}[Prolongation]
Let function $f(t,\,u)$ is defined and continuous in domain $\D\,$. If a solution $u\ =\ \phi(t)$ to problem \eqref{eq:ODE} exists in the interval $t\ \in\ [t_0,\,\alpha)$ cannot be prolongated beyond the point $t\ =\ \alpha$, then it can happen only for one of three following reasons:
\begin{enumerate}
  \item $\alpha\ =\ +\infty\,$,
  \item When $t\ \to\ \alpha\ -\ 0$, $\abs{\phi(t)}\ \to\ +\infty\,$,
  \item When $t\ \to\ \alpha\ -\ 0$, distance between the point $\bigl(t,\,\phi(t)\bigr)$ to the boundary $\partial\D$ goes to zero.
\end{enumerate}
\end{theorem}

\begin{theorem}[Uniqueness, \cite{Osgood1898}]
If function $f(t,\,u)$ is defined in domain $\D$ and for any pair of points $(t_1,\,u_1)\,$, $(t_2,\,u_2)\ \in\ \D$ satisfies the condition
\begin{equation}\label{eq:lipschitz}
  \abs{f(t_1,\,u_1)\ -\ f(t_2,\,u_2)}\ \leq\ \omega\,\bigl(\abs{u_1\ -\ u_2}\bigr)\,,
\end{equation}
where $\omega(u)\ >\ 0$ is continuous for $0\ <\ u\ \leq\ \U$ and
\begin{equation*}
  \int_{t_0\ +\ \eps}^{\,\U}\,\frac{\ud u}{\omega(u)}\ \to\ +\infty\,, \qquad \mbox{ as }\ \qquad \eps\ \to\ 0\,,
\end{equation*}
then through any point $(t_0,\,u_0)\ \in\ \D$ passes at most one integral curve of equation \eqref{eq:ODE}.
\end{theorem}
Suitable functions $\omega(u)$ which satisfy the conditions of the Theorem are
\begin{align*}
  \omega(u)\ &=\ \K\,u\,, \\
  \omega(u)\ &=\ \K\,u\,\abs{\ln u}\,, \\
  \omega(u)\ &=\ \K\,u\,\abs{\ln u}\cdot\ln\abs{\ln u}\,, \\
  \omega(u)\ &=\ \K\,u\,\abs{\ln u}\cdot\ln\abs{\ln u}\cdot\ln\ln\abs{\ln u}\,, \qquad \mbox{ \etc}
\end{align*}
Above $\K$ is a positive constant. If we take $\omega(u)\ =\ \K\,u$ then condition \eqref{eq:lipschitz} becomes the well-known \textsc{Lipschitz}\footnote{Rudolf Otto Sigismund \textsc{Lipschitz} (1832 -- 1903) is a German Mathematician who made contributions to Mathematical Analysis. He studied with Gustav \textsc{Dirichlet} at the University of Berlin.} condition for function $f(t,\,u)$ in the second variable $u\,$. In order to satisfy the \textsc{Lipschitz} condition, it is sufficient for function $f(t,\,u)$ to be defined in a domain $\D$ convex in $u$ and to have a bounded derivative $\pd{f}{u}$ in this domain. Later, \textsc{Wintner} showed that \textsc{Osgood}\footnote{William Fogg \textsc{Osgood} (1864 -- 1943) was an American Mathematician born un Boston, MA. He studied in Universities of G\"ottingen and Erlangen, Germany. He made contributions to Mathematical and Complex Analysis, Conformal Mappings.}'s Theorem conditions are sufficient for the convergence of \textsc{Picard}\footnote{Charles \'Emile \textsc{Picard} (1856 -- 1941) was a French Mathematician who made contributions to the Mathematical Analysis. He was married to Marie, a daughter of his Professor Charles \textsc{Hermite}.}'s iterations to a local solution on a sufficiently small interval \cite{Wintner1946}.

\subsubsection{Counterexamples}

The Theorem above gives the existence of solutions only \emph{locally} in time $\bigl[\,t_0,\,t_0\ +\ \delta\,\bigr)$. The only reason for a solution not to exist beyond the given interval is the \emph{blow-up} (\ie~the solution becomes unbounded). Otherwise, the solution exists globally. For instance, the following problem
\begin{equation*}
  \dot{u}\ =\ 1\ +\ u^2\,, \qquad u(0)\ =\ 0\,,
\end{equation*}
has the unique exact solution $u(t)\ =\ \tan(t)$ which exists only in the interval $\bigl[\,0,\,\frac{\pi}{2}\bigr)\,$. At time $t\ =\ \frac{\pi}{2}$ the solution blows up.

The uniqueness property is sometimes violated as well. For instance, the scalar equation
\begin{equation*}
  \dot{u}\ =\ \sqrt{\abs{u}}\,, \qquad u(0)\ =\ 0\,,
\end{equation*}
has actually infinitely many solutions. Two examples are $u(t)\ \equiv\ 0$ and $u(t)\ =\ \dfrac{t^2}{4}\,$. Obviously, the right-hand side does not satisfy the \textsc{Lipschitz} condition which guarantees the uniqueness.


\subsection*{Acknowledgments}
\addcontentsline{toc}{subsection}{Acknowledgments}

The Author would like to thank Professor Didier~\textsc{Clamond} (University of Nice Sophia Antipolis, France) for introducing me to \textsc{Fourier}-type pseudo-spectral methods. I would like also to thank Professor Laurent~\textsc{Gosse} (CNR, Italy) who brought Author's attention to \textsc{Trefftz} methods. Special thanks go to my colleague Prof.~Paul-\'Eric \textsc{Chaudru De Raynal} (LAMA, Universit\'e Savoie Mont Blanc, France) for his help with Monte--Carlo simulations. Finally, the Author would like to thank Professors Nathan~\textsc{Mendes} and Marx~\textsc{Chhay} for giving me an opportunity to deliver these lectures in Brazil.


\addcontentsline{toc}{section}{References}
\bibliographystyle{abbrv}
\bibliography{biblio}

\begin{thebibliography}{10}

\bibitem{Benjamin2010}
A.~T. Benjamin and D.~Walton.
\newblock {Combinatorially composing Chebyshev polynomials}.
\newblock {\em Journal of Statistical Planning and Inference},
  140(8):2161--2167, aug 2010.

\bibitem{Boyd2000}
J.~P. Boyd.
\newblock {\em {Chebyshev and Fourier Spectral Methods}}.
\newblock New York, 2nd edition, 2000.

\bibitem{Caflisch1998}
R.~E. Caflisch.
\newblock {Monte Carlo and quasi-Monte Carlo methods}.
\newblock {\em Acta Numerica}, 7:1--49, 1998.

\bibitem{Cheung1989}
Y.~K. Cheung, W.~G. Jin, and O.~C. Zienkiewicz.
\newblock {Direct solution procedure for solution of harmonic problems using
  complete, non-singular, Trefftz functions}.
\newblock {\em Comm. Appl. Num. Meth.}, 5:159--169, 1989.

\bibitem{Chhay2015}
M.~Chhay, D.~Dutykh, M.~Gisclon, and C.~Ruyer-Quil.
\newblock {Asymptotic heat transfer model in thin liquid films}.
\newblock {\em Submitted}, pages 1--24, 2015.

\bibitem{Cooley1965}
J.~W. Cooley and J.~W. Tukey.
\newblock {An algorithm for the machine calculation of complex Fourier series}.
\newblock {\em Mathematics of Computation}, 19(90):297--297, may 1965.

\bibitem{Dormand1980}
J.~R. Dormand and P.~J. Prince.
\newblock {A family of embedded Runge-Kutta formulae}.
\newblock {\em J. Comp. Appl. Math.}, 6:19--26, 1980.

\bibitem{Einstein1905}
A.~Einstein.
\newblock {{\"{U}}ber die von der molekularkinetischen Theorie der W{\"{a}}rme
  geforderte Bewegung von in ruhenden Fl{\"{u}}ssigkeiten suspendierten
  Teilchen}.
\newblock {\em Annalen der Physik}, 322(8):549--560, 1905.

\bibitem{Einstein1926}
A.~Einstein.
\newblock {Die Ursache der M{\"{a}}anderbildung der Flu{\ss}l{\"{a}}ufe und des
  sogenannten Baerschen Gesetzes}.
\newblock {\em Die Naturwissenschaften}, 14(11):223--224, 1926.

\bibitem{Evans2010}
L.~C. Evans.
\newblock {\em {Partial Differential Equations}}.
\newblock American Mathematical Society, Providence, Rhode Island, 2 edition,
  2010.

\bibitem{Fornberg1996}
B.~Fornberg.
\newblock {\em {A practical guide to pseudospectral methods}}.
\newblock Cambridge University Press, Cambridge, 1996.

\bibitem{Fourier1822}
J.~Fourier.
\newblock {\em {Th{\'{e}}orie analytique de la chaleur}}.
\newblock Didot, Paris, 1822.

\bibitem{Gisclon1998}
M.~Gisclon.
\newblock {A propos de l'{\'{e}}quation de la chaleur et de l'analyse de
  Fourier}.
\newblock {\em Le journal de maths des {\'{e}}l{\`{e}}ves}, 1(4):190--197,
  1998.

\bibitem{Griffiths2015}
D.~F. Griffiths, J.~W. Dold, and D.~J. Silvester.
\newblock {Separation of Variables}.
\newblock In {\em Essential Partial Differential Equations}, pages 129--159.
  Springer International Publishing, 2015.

\bibitem{Hadamard1902}
J.~Hadamard.
\newblock {Sur les probl{\`{e}}mes aux d{\'{e}}riv{\'{e}}es partielles et leur
  signification physique}.
\newblock {\em Princeton University Bulletin}, pages 49--52, 1902.

\bibitem{Hairer2002}
E.~Hairer, C.~Lubich, and G.~Wanner.
\newblock {\em {Geometric Numerical Integration}}, volume~31 of {\em Spring
  Series in Computational Mathematics}.
\newblock Springer-Verlag, Berlin, Heidelberg, second edition, 2006.

\bibitem{Hairer2009}
E.~Hairer, S.~P. N{\o}rsett, and G.~Wanner.
\newblock {\em {Solving ordinary differential equations: Nonstiff problems}}.
\newblock Springer, 2009.

\bibitem{Hairer1996}
E.~Hairer and G.~Wanner.
\newblock {\em {Solving Ordinary Differential Equations II. Stiff and
  Differential-Algebraic Problems}}.
\newblock Springer Series in Computational Mathematics, Vol. 14, 1996.

\bibitem{Herrera1982}
I.~Herrera.
\newblock {Boundary methods: development of complete systems of solutions}.
\newblock In T.~Kawai, editor, {\em Finite Elements Flow Analysis}, pages
  897--906, Tokyo, 1982. University of Tokyo Press.

\bibitem{Herrera1984}
I.~Herrera.
\newblock {\em {Boundary Methods: An Algebraic Theory}}.
\newblock Pitman, 1984.

\bibitem{Higham2001}
D.~J. Higham.
\newblock {An Algorithmic Introduction to Numerical Simulation of Stochastic
  Differential Equations}.
\newblock {\em SIAM Review}, 43(3):525--546, jan 2001.

\bibitem{Jaynes2003}
E.~T. Jaynes.
\newblock {\em {Probability Theory}}.
\newblock Cambridge University Press, Cambridge, 2003.

\bibitem{Kita1995}
E.~Kita and N.~Kamiya.
\newblock {Trefftz method: an overview}.
\newblock {\em Advances in Engineering Software}, 24:3--12, 1995.

\bibitem{Kreiss1992}
H.~O. Kreiss and G.~Scherer.
\newblock {Method of lines for hyperbolic equations}.
\newblock {\em SIAM Journal on Numerical Analysis}, 29:640--646, 1992.

\bibitem{Lapeyre2011}
B.~Lapeyre and J.~Lelong.
\newblock {A framework for adaptive Monte-Carlo procedures}.
\newblock {\em Monte Carlo Methods Appl.}, 17(1):77--98, 2011.

\bibitem{Mendes2005}
N.~Mendes and P.~C. Philippi.
\newblock {A method for predicting heat and moisture transfer through
  multilayered walls based on temperature and moisture content gradients}.
\newblock {\em Int. J. Heat Mass Transfer}, 48(1):37--51, 2005.

\bibitem{Osgood1898}
W.~F. Osgood.
\newblock {Beweis der Existenz einer L{\"{o}}sung der Differentialgleichung
  $\frac{\ud y}{\ud x}\ =\ f(x,\,y)$ ohne Hinzunahme der
  Cauchy-Lipschitz'schen Bedingung}.
\newblock {\em Monatshefte f{\"{u}}r Mathematik und Physik}, 9(1):331--345, dec
  1898.

\bibitem{Pearson1905}
K.~Pearson.
\newblock {The Problem of the Random Walk}.
\newblock {\em Nature}, 72:294, 1905.

\bibitem{Philibert2005}
J.~Philibert.
\newblock {One and a half century of diffusion: Fick, Einstein, before and
  beyond}.
\newblock In J.~K{\"{a}}rger, F.~Grinberg, and P.~Heitjans, editors, {\em
  Diffusion Fundamentals}, pages 8--17. Leipzig Universt{\"{a}}tsverlag,
  Leipzig, 2005.

\bibitem{Reddy1992}
S.~C. Reddy and L.~N. Trefethen.
\newblock {Stability of the method of lines}.
\newblock {\em Numerische Mathematik}, 62(1):235--267, 1992.

\bibitem{Schiesser1994}
W.~E. Schiesser.
\newblock {Method of lines solution of the Korteweg-de vries equation}.
\newblock {\em Computers Mathematics with Applications}, 28(10-12):147--154,
  1994.

\bibitem{Shampine1994}
L.~F. Shampine.
\newblock {ODE solvers and the method of lines}.
\newblock {\em Numerical Methods for Partial Differential Equations},
  10(6):739--755, 1994.

\bibitem{Shampine1997}
L.~F. Shampine and M.~W. Reichelt.
\newblock {The MATLAB ODE Suite}.
\newblock {\em SIAM Journal on Scientific Computing}, 18:1--22, 1997.

\bibitem{Solin2005}
P.~Solin.
\newblock {\em {Partial Differential Equations and the Finite Element Method}}.
\newblock John Wiley {\&} Sons, Inc., Hoboken, New Jersey, 2005.

\bibitem{Trefethen2000}
L.~N. Trefethen.
\newblock {\em {Spectral methods in MatLab}}.
\newblock Society for Industrial and Applied Mathematics, Philadelphia, PA,
  USA, 2000.

\bibitem{Trefftz1926}
E.~Trefftz.
\newblock {Gegenst{\"{u}}ck zum ritzschen Verfahren}.
\newblock In {\em Proc. 2nd Int. Cong. Appl. Mech.}, pages 131--137,
  Z{\"{u}}rich, 1926.

\bibitem{Uecker2009}
H.~Uecker.
\newblock {A short ad hoc introduction to spectral methods for parabolic PDE
  and the Navier-Stokes equations}.
\newblock Technical report, Carl von Ossietzky Universit{\"{a}}t Oldenburg,
  Oldenburg, Germany, 2009.

\bibitem{Leer2006}
B.~van Leer.
\newblock {Upwind and High-Resolution Methods for Compressible Flow: From Donor
  Cell to Residual-Distribution Schemes}.
\newblock {\em Commun. Comput. Phys.}, 1:192--206, 2006.

\bibitem{Vertesi1990}
P.~V{\'{e}}rtesi.
\newblock {Optimal Lebesgue constant for Lagrange interpolation}.
\newblock {\em SIAM J. Numer. Anal.}, 27:1322--1331, 1990.

\bibitem{Wintner1946}
A.~Wintner.
\newblock {On the Convergence of Successive Approximations}.
\newblock {\em American Journal of Mathematics}, 68(1):13, jan 1946.

\end{thebibliography}

\end{document}